\documentclass[11pt,a4paper]{amsart} 
\usepackage{amssymb}
\usepackage[all]{xy}
\usepackage{subfig}
\SelectTips{cm}{}
\calclayout 

\title[Representations of quivers]{Representations of quivers
via reflection functors}
\author[Henning Krause]{Henning Krause}
\address{Fakult\"at f\"ur Mathematik\\
Universit\"at Bielefeld\\ 33501 Bielefeld\\ Germany}
\email{hkrause@math.uni-bielefeld.de}

\newtheorem{lem}{Lemma}[subsection]
\newtheorem{prop}[lem]{Proposition}
\newtheorem{cor}[lem]{Corollary}
\newtheorem{thm}[lem]{Theorem}

\newtheorem{Lem}{Lemma}[section] 

\theoremstyle{definition}
\newtheorem{exm}[lem]{Example}
\newtheorem{defn}[lem]{Definition}

\newtheorem{rem}[lem]{Remark}

\numberwithin{equation}{subsection}

\newcommand{\smatrix}[1]{\left[\begin{smallmatrix}#1\end{smallmatrix}\right]}

\renewcommand{\leq}{\leqslant}
\renewcommand{\geq}{\geqslant}
\newcommand{\Aut}{\operatorname{Aut}\nolimits}
\newcommand{\Rep}{\operatorname{Rep}\nolimits}
\newcommand{\rep}{\operatorname{rep}\nolimits}

\newcommand{\card}{\operatorname{card}\nolimits}
\newcommand{\Irr}{\operatorname{Irr}\nolimits}
\newcommand{\Rad}{\operatorname{Rad}\nolimits}
\newcommand{\rad}{\operatorname{rad}\nolimits}

\newcommand{\id}{\operatorname{id}\nolimits}

\newcommand{\End}{\operatorname{End}\nolimits}
\newcommand{\Hom}{\operatorname{Hom}\nolimits}

\renewcommand{\Im}{\operatorname{Im}\nolimits}
\newcommand{\Ker}{\operatorname{Ker}\nolimits}
\newcommand{\Coker}{\operatorname{Coker}\nolimits}

\renewcommand{\dim}{\operatorname{dim}\nolimits}

\newcommand{\Ext}{\operatorname{Ext}\nolimits}

\newcommand{\op}{\mathrm{op}}
\newcommand{\Id}{\mathrm{Id}}

\newcommand{\lto}{\longrightarrow}
\newcommand{\xto}{\xrightarrow}

\def\a{\alpha}
\def\b{\beta}
\def\e{\varepsilon}
\def\d{\delta}
\def\g{\gamma}
\def\p{\phi}
\def\r{\rho}
\def\s{\sigma}
\def\t{\tau}

\def\m{\mu}

\def\la{\lambda}
\def\De{\varDelta}
\def\Ga{\varGamma}
\def\La{\varLambda}

\def\bbN{\mathbb N}
\def\bbP{\mathbb P}
\def\bbQ{\mathbb Q}

\def\bbZ{\mathbb Z}

\begin{document}
\maketitle 

\setcounter{tocdepth}{1}
\tableofcontents

These are the notes for a course on representations of quivers for
second year students in Paderborn in summer 2007. My aim was to
provide a basic introduction without using any advanced methods. It
turns out that a good knowledge of linear algebra is sufficient for
proving Gabriel's theorem.  Thus we classify the quivers of finite
representation type and study their representations using reflection
functors. The course was complemented by problem sessions run by
Andrew Hubery; see his homepage for interesting exercises and further
material. I wish to thank him for many useful discussions on the
subject and the students of this course for their enthusiasm.

Further material has been added after giving this course in Bielefeld
in summer 2010. This includes a discussion of regular representations
and wild phenomena. In particular, two classical examples are covered:
representations of the Kronecker quiver and representations of the
Klein four group.

\section{Representations of quivers}

In this section we introduce our basic concepts: quivers and their
representations. Throughout we fix a field $k$.

\subsection{Quivers}
A \emph{quiver} is a directed graph, which is assumed to be finite.
More precisely, a quiver is a quadruple $Q=(Q_0,Q_1,s,t)$ consisting
of a finite set $Q_0$ of \emph{vertices}, a finite set $Q_1$ of
\emph{arrows}, and two maps $s,t\colon Q_1\to Q_0$. An arrow $\a\in
Q_1$ \emph{starts} at $s(\a)$ and \emph{terminates} at $t(\a)$.  We
sometimes write $\a\colon s(\a)\to t(\a)$.

A \emph{non-trivial path} of length $r\ge 1$ in $Q$ is a sequence
$\xi=\xi_r\ldots\xi_1$ of arrows satisfying $t(\xi_{p})=s(\xi_{p+1})$
for $1\le p < r$.  We write
$$i_1\stackrel{\xi_1}\lto i_2\stackrel{\xi_2}\lto \cdots
\stackrel{\xi_r}\lto i_{r+1}.$$ The path
$\xi$ starts at $s(\xi)=s(\xi_1)$ and terminates at
$t(\xi)=t(\xi_r)$. For each vertex $i$, we have in addition the
\emph{trivial path} $\e_i$ of length zero with $s(\e_i)=i=t(\e_i)$.

For a pair $i,j$ of vertices, let $Q(i,j)$ denote the set of paths
$\xi$ with $s(\xi)=i$ and $t(\xi)=j$. The obvious composition of paths
induces maps
$$Q(i,\xi)\colon Q(i,s(\xi))\lto Q(i,t(\xi))\quad\text{and}\quad
Q(\xi,j)\colon Q(t(\xi),j)\lto Q(s(\xi),j)$$
with $Q(i,\xi)(\mu)=\xi\mu$ and $Q(\xi,j)(\nu)=\nu\xi$.

\subsection{Representations}
Let $Q$ be a quiver. A \emph{representation} of $Q$ is a collection
$$X=(X_i,X_\a)_{i\in Q_0,\a\in Q_1}$$ consisting of a vector space $X_i$
for each vertex $i$ and a linear map $X_\a\colon X_{s(\a)}\to
X_{t(\a)}$ for each arrow $\a$. A \emph{morphism} of representations
$\p\colon X\to Y$ is a collection $\p=(\p_i)_{i\in Q_0}$ of linear
maps $\p_i\colon X_i\to Y_i$ for each vertex $i$ such that
$Y_\a\p_{s(\a)}=\p_{t(\a)}X_\a$ for each arrow $\a$. In other words,
for each arrow $\a$ we have a commutative diagram
$$\xymatrix{
X_{s(\a)}\ar[d]^{X_\a}\ar[rr]^{\p_{s(\a)}}&&Y_{s(\a)}\ar[d]^{Y_\a}\\
X_{t(\a)}\ar[rr]^{\p_{t(\a)}}&&Y_{t(\a)} }$$ The \emph{composition} of
$\p$ with $\psi\colon Y\to Z$ is given by $(\psi\p)_i=\psi_i\p_i$ for
each vertex $i$. For each representation $X$, we have the
\emph{identity morphism} $\id_X\colon X\to X$ with $(\id_X)_i=\id_{X_i}$
for all $i$. The set of morphisms $X\to Y$ we denote by $\Hom(X,Y)$
and we write $\End(X)$ for the set of endomorphisms $X\to X$.

This defines a category $\Rep (Q,k)$. We denote by $\rep (Q,k)$ the
full subcategory with objects the finite dimensional
representations. Here, a representation $X$ is \emph{finite
dimensional} if each vector space $X_i$ is finite dimensional.

\subsection{The category of representations}
Various concepts that are defined for vector spaces carry over to
representations of $Q$ by applying the vector space definition
\emph{point-wise}, that is, for each vertex $i$. 

Fix a pair $X,Y$ of representations of $Q$.  We call $X$ a
\emph{subrepresentation} of $Y$ and write $X\subseteq Y$, if $X_i$ is
a subspace of $Y_i$ for each vertex $i$ and $X_\a(x)=Y_\a(x)$ for each
arrow $\a$ and $x\in X_{s(\a)}$.

Given a morphism $\p\colon X\to Y$, its \emph{kernel} $\Ker\p$ is by
definition the subrepresentation of $X$ with $(\Ker\p)_i=\Ker\p_i$ for
each vertex $i$. The \emph{cokernel} $\Coker\p$ and the \emph{image}
$\Im\p$ are defined analogously. Note that $\p$ is a
\emph{monomorphism} if and only if $\Ker\p=0$, while $\p$ is an
\emph{epimorphism} if and only if $\Coker\p=0$. The morphism $\p$ is
an \emph{isomorphism} if each $\p_i$ is an isomorphism. One defines
addition and scalar multiplication for morphisms $X\to Y$ point-wise
and that makes $\Hom(X,Y)$ into a vector space.

The \emph{dimension vector} of a finite dimensional representation $X$
is the vector $\dim X$ in $\bbZ^{Q_0}$ with
$$(\dim X)_i=\dim X_i\quad (i\in Q_0).$$

\subsection{Some special quivers}
We consider three quivers in more detail. 
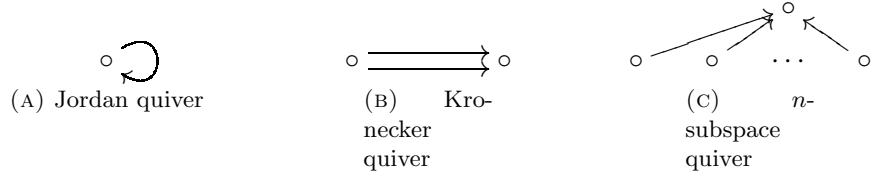
\begin{figure}[!ht]
\centering 
\subfloat[Jordan quiver]
{\makebox[11em]{$\xymatrix@!0 @R=1.8em @C=2.6em {{}\\ \circ\ar@(ur,dr)}$}}
\subfloat[Kronecker quiver]
{\makebox[11em]{$\xymatrix@!0 @R=1.8em @C=2.6em {{}\\ \circ\ar@<3.0pt>[rr]\ar@<-3.0pt>[rr]&&\circ}$}}
\subfloat[$n$-subspace quiver]
{\makebox[11em]{$\xymatrix@!0 @R=1.8em @C=2.6em {
&&\circ\\ \circ\ar[urr]&\circ\ar[ur]&\cdots&\circ\ar[ul]}$}}
\end{figure}
The representations of the \emph{Jordan quiver} are endomorphisms of a
single vector space and their classification can be formulated in
terms of Jordan blocks. The representations of the \emph{Kronecker
quiver} are pairs of linear maps up to simultaneous conjugation.  The
\emph{n-subspace quiver} has $n+1$ vertices and its representations
are basically configurations of $n$ subspaces of a fixed vector space.

\subsection{Duality}
Let $Q^\op$ denote the \emph{opposite quiver} which is obtained from
$Q$ by reversing all arrows. The vector space duality $D=\Hom_k(-,k)$
induces a \emph{duality}
$$D\colon \Rep(Q,k)\lto\Rep(Q^\op,k).$$ Given a representation $X$ of
$Q$, let $(DX)_i=D(X_i)$ and $(DX)_\a=D(X_\a)$ for $i\in Q_0$ and
$\a\in Q_1$. For a morphism $\p\colon X\to Y$, let
$(D\p)_i=D(\p_i)$.

Let $V,W$ be a pair of vector spaces. Recall that there is a canonical
monomorphism $\e_V\colon V\to D^2V$ defined by $\e_V(x)(\p)=\p(x)$ for
$x\in V$ and $\p\in DV$. This induces an isomorphism
$$\Hom(W,DV)\xto{\sim}\Hom(V,DW)$$
by sending $\p\colon W\to DV$ to $(D\p)\e_V$.
\begin{lem}
\label{le:dual}
For $X\in\Rep(Q,k)$ and $Y\in\Rep(Q^\op,k)$, there is a canonical
monomorphism $\e_X\colon X\to D^2X$ which induces a natural
isomorphism
$$\Hom(Y,DX)\xto{\sim}\Hom(X,DY)$$ by sending $\p\colon Y\to DX$ to
$(D\p)\e_X$.
\end{lem}
\begin{proof}
Use the linear maps $\e_{X_i}$ and
$\Hom(Y_i,DX_i)\xto{\sim}\Hom(X_i,DY_i)$.
\end{proof}

\subsection{Simple representations}
A \emph{simple} representation is defined to be a non-zero
representation with no proper subrepresentations.

Given a vertex $i$, let $S(i)$ be the representation with 
$$S(i)_j=\begin{cases} k&\text{if }j=i,\\0&\text{if }j\neq i,
\end{cases}\quad\text{and}\quad S(i)_\a=0$$
for $j\in Q_0$ and $\a\in Q_1$. This representation is simple. 

\begin{lem}\label{le:simple}
Let $X$ be a representation and suppose that $i$ is a vertex
with $X_i\neq 0$ and $X_\a=0$ for each arrow $\a$ starting at
$i$. Then $S(i)$ is a subrepresentation of $X$.
\end{lem}
\begin{proof}
The assumption on $X$ implies  $\Hom(S(i),X)\cong X_i$.
\end{proof}

Suppose that $Q$ has no \emph{oriented cycles}, that is, non-trivial
paths from a vertex to itself. Then for any simple representation $S$
of $Q$, there exists a unique vertex $i$ such that $S\cong S(i)$. This
follows from Lemma~\ref{le:simple}. On the other hand, there are
additional simple representations if $Q$ has oriented cycles.

\begin{exm}
The finite dimensional simple representations of the Jordan quiver are
parametrised by the monic irreducible polynomials over $k$. More
precisely, the representation corresponding to such a polynomial
$\sum_{i=0}^d\la_i t^i$ of degree $d$ is the pair $(X,\p)$ consisting of
the vector space $X=k^d$ and the endomorphism $\p\colon X\to X$ with
$\p(e_i)=e_{i+1}$ for $1\le i <d$ and $\p(e_d)=\sum_{i=1}^d-\la_{i-1}e_i$.
\end{exm}

\subsection{Projective and injective representations}
For each vertex $i$, we define a \emph{projective representation}
$P(i)$ and an \emph{injective representation} $I(i)$.

Given any set $X$, we denote by $k[X]$ the vector space with basis $X$,
that is, $k[X]$ is the set of linear combinations $\sum_p\lambda_px_p$
with $x_p\in X$, $\lambda_p\in k$, and almost all $\lambda_p=0$. For a map
$\p\colon X\to Y$, let $k[\p]\colon k[X]\to k[Y]$ be the linear map sending
$\sum_p\lambda_px_p$ to $\sum_p\lambda_p\p(x_p)$. 

Now define
$$P(i)_j=k[Q(i,j)]\quad\text{and}\quad P(i)_\a=k[Q(i,\a)]$$
for $j\in Q_0$ and $\a\in Q_1$. Dually, we define
$$I(i)_j=Dk[Q(j,i)]\quad\text{and}\quad I(i)_\a=Dk[Q(\a,i)].$$ Note that
$I(i)=D\bar P(i)$ where $\bar P(i)$ refers to the projective representation of
$Q^\op$.

\begin{lem}
\label{le:hom}
Let $X$ be a representation of $Q$. Then there are natural isomorphisms
$$\Hom(P(i),X)\cong X_i\quad\text{and}\quad
\Hom(X,I(i))\cong DX_i.$$
\end{lem}
\begin{proof}
The isomorphism $\Hom(P(i),X)\to X_i$ sends $\p$ to $\p_i(\e_i)$. Its inverse
map sends $x\in X_i$ to the morphism $\p\colon P(i)\to X$ with
$$\p_j(\xi)=X_{\xi_r}\ldots X_{\xi_1}(x)$$ for a basis element
$\xi=\xi_r\ldots\xi_1$ of $P(i)_j$.

The isomorphism $\Hom(X,I(i))\cong DX_i$ follows from the first using
Lemma~\ref{le:dual} since
\begin{equation*}
\Hom(X,I(i))=\Hom(X,D\bar P(i))\cong\Hom(\bar P(i),DX)\cong DX_i.\qedhere
\end{equation*}
\end{proof}

\begin{lem}
\label{le:proj}
\begin{enumerate} 
\item The representations $P(i)$ ($i\in Q_0$) are pairwise non-isomorphic.
\item The representations $I(i)$ ($i\in Q_0$) are pairwise non-isomorphic.
\end{enumerate}
\end{lem}
\begin{proof}
Fix two vertices $i,j$ and suppose $P(j)\cong P(i)$. Then
$$k=S(i)_i\cong\Hom(P(i),S(i))\cong\Hom(P(j),S(i))\cong
S(i)_j=\begin{cases} k&\text{if }j=i,\\0&\text{if }j\neq i,
\end{cases}$$
by Lemma~\ref{le:hom}. Thus $j=i$. The proof for the $I(i)$ is analogous.
\end{proof}

\begin{lem}
\label{le:endproj}
Suppose $Q$ has no oriented cycles. Then $P(i)$ and $I(i)$ are finite
dimensional with $\End(P(i))\cong k\cong \End(I(i))$.
\end{lem}
\begin{proof}
If $Q$ has no oriented cycles then $Q(i,j)$ is finite for all pairs of
vertices $i,j$ because the quiver is finite. We have $$\End(P(i))\cong P(i)_i=k[Q(i,i)]=k$$
by Lemma~\ref{le:hom}, and a similar argument works for $I(i)$.
\end{proof}

\begin{rem} 
\label{re:inj}
The representation $P(i)$ is a \emph{projective object} in the sense that
for every epimorphism $X\to Y$ the induced map $\Hom(P(i),X)\to\Hom(P(i),Y)$
is surjective. Dually, $I(i)$ is an \emph{injective object} in the sense that
for every monomorphism $X\to Y$ the induced map $\Hom(Y,I(i))\to\Hom(X,I(i))$
is surjective. This follows from Lemma~\ref{le:hom}.
\end{rem}

\section{Direct sum decompositions}
\label{se:krs}
In this section we consider finite dimensional representations. We
show that each representation decomposes essentially uniquely into
indecomposable representations.

\subsection{Direct sums}
Let $X_1,\ldots,X_r$ be a finite number of representations. A
\emph{direct sum}
$$X=X_1\oplus\ldots\oplus X_r$$ is a representation $X$ together with
morphisms $\iota_i\colon X_i\to X$ and $\pi_i\colon X\to X_i$ for
$1\leq i\leq r$ such that $\sum_{i=1}^r\iota_i\pi_i=\id_X$ and
$\pi_i\iota_i=\id_{X_i}$ for all $i$. Note that we can identify each
$X_i$ via $\iota_i$ with a subrepresentation of $X$. Then we obtain
\begin{equation}
\label{eq:sum}
X=\sum_{i=1}^rX_i\quad\text{and}\quad X_i\cap\sum_{i'\neq i} X_{i'}=0
\quad\text{for}\quad 1\le i\le r.
\end{equation} Here $\sum_{i\in I}X_i$ refers to the smallest
subrepresentation of $X$ containing $X_i$ for all $i\in I$.
Conversely, if $X_1,\ldots,X_r$ is a family of subrepresentations of a
representation $X$ satisfying \eqref{eq:sum} then
$X=X_1\oplus\ldots\oplus X_r$. In that case we take for $\iota_i\colon
X_i\to X$ the inclusion morphism and let
$\pi_i=(\r_i\iota_i)^{-1}\r_i$, where $\r_i\colon X\to X/{\sum_{i'\neq
i} X_{i'}}$ denotes the canonical morphism.

A family of subrepresentations $X_1,\ldots,X_r$ of $X$ satisfying
\eqref{eq:sum} is called a \emph{direct sum decomposition} of $X$.
Note that one can check \eqref{eq:sum} point-wise for each vertex.
Thus we have a decomposition $X=X_1\oplus\ldots\oplus X_r$ if and only
if we have a vector space decomposition $X_j=(X_1)_j\oplus\ldots\oplus
(X_r)_j$ for each vertex $j$.

\begin{lem} 
\label{le:sum}
Let $X=X_1\oplus\ldots\oplus X_r$ and $Y=Y_1\oplus\ldots\oplus Y_s$.
Then we have induced vector space decompositions
$$\bigoplus_{i=1}^r\Hom(X_i,Y)=\Hom(X,Y)=\bigoplus_{j=1}^s\Hom(X,Y_j).$$
\end{lem}
\begin{proof}
Let $\iota^*_i=\Hom(\iota_i,Y)$ and $\pi^*_i=\Hom(\pi_i,Y)$ for $1\leq
i\leq r$. Then we have $\sum_{i=1}^r\pi^*_i\iota^*_i=\id_{\Hom(X,Y)}$
and $\iota^*_i\pi^*_i=\id_{\Hom(X_i,Y)}$ for all $i$. This proves the
first equality. The argument for the second equality is analogous.
\end{proof}

It follows from Lemma~\ref{le:sum} that a direct sum of
$X_1,\ldots,X_r$ is unique up to an isomorphism.  Thus we may speak of
\emph{the} direct sum and the notation $X_1\oplus\ldots\oplus X_r$ is
well-defined. We write $X^r=X\oplus\ldots\oplus X$ for the direct sum of $r$
copies of a representation $X$.

A representation $X$ is called \emph{indecomposable} if $X\neq 0$ and
$X=X_1\oplus X_2$ implies $X_1=0$ or $X_2=0$.

\subsection{Fitting's lemma}
We fix a representation $X$ and study the ring of endomorphisms $\End(X)$.

\begin{lem} 
\label{le:fitting}
Let $X$ be a representation and $\p$ an endomorphism.
\begin{enumerate}
\item For large enough $r$, we have  $X=\Im\p^r\oplus\Ker\p^r$.
\item If $X$ is indecomposable, then $\p$ is either an automorphism or
nilpotent.
\end{enumerate}
\end{lem}
\begin{proof}
  Because $X$ is finite dimensional, we may choose $r$ large enough so
  that $\Im \p^r=\Im\p^{r+1}$. Thus $\p^r\colon\Im\p^r\to\Im\p^{2r}$
  is an isomorphism and we denote by $\psi$ its inverse.  Furthermore, let
  $\iota_1\colon\Im\p^r\to X$ and $\iota_2\colon \Ker\p^r\to X$ denote
  the inclusions. We put $\pi_1=\psi\p^r\colon X\to \Im\p^r$ and
  $\pi_2=\id_X-\psi\p^r\colon X\to \Ker\p^r$. Then
  $\iota_1\pi_1+\iota_2\pi_2=\id_X$ and $\pi_i\iota_i=\id_{X_i}$ for
  $i=1,2$. Thus $X=\Im\p^r\oplus\Ker\p^r$.  Part (2) is an immediate
  consequence of (1).
\end{proof}

A ring is called \emph{local} if the sum of two non-units is
again a non-unit.

\begin{prop}
\label{pr:fitting}
A representation $X$ is indecomposable if and only if $\End(X)$ is local.
\end{prop}
\begin{proof}
Let $X$ be indecomposable and $\p,\p'\in\End(X)$.  Suppose $\p+\p'$ is
invertible, say $\r(\p+\p')=\id_X$. If $\p$ is non-invertible then
$\r\p$ is non-invertible. Thus $\r\p$ is nilpotent, say $(\r\p)^r=0$,
by Lemma~\ref{le:fitting}. We obtain
$$(\id_X-\r\p)(\id_X+\r\p+\ldots+(\r\p)^{r-1})=\id_X.$$ Therefore
$\r\p'=\id_X-\r\p$ is invertible whence $\p'$ is invertible.

If $X=X_1\oplus X_2$ with $X_i\neq 0$ for $i=1,2$, then we have
idempotent endomorphisms $\e_i$ of $X$ with $\Im\e_i=X_i$. Clearly,
each $\e_i$ is non-invertible but $\id_X=\e_1+\e_2$.
\end{proof}

The assumption on $X$ to be finite dimensional is necessary.
\begin{exm}
Let $k[t]$ denote the polynomial ring in one variable and consider the
following representation of the Kronecker quiver. $$\xymatrix@!0
@R=1.8em @C=2.6em {X:&k[t]\ar@<3.0pt>[rr]^-{\cdot
t}\ar@<-3.0pt>[rr]_-{\id}&&k[t]}$$ The endomorphism ring of $X$ is
isomorphic to $k[t]$. Thus $X$ is indecomposable but $\End(X)$ is not
local.
\end{exm}

\subsection{The radical}
\label{se:radical}
Given a pair $X,Y$ of representations, we define the \emph{radical}
$$\Rad(X,Y)=\left\{\p\in\Hom(X,Y)\left|
\genfrac{}{}{0pt}{}{\t\p\s\text{ is non-invertible for every pair
}Z\xto{\s}X}{\text{and }Y\xto{\t}Z\text{ with }Z\text{
indecomposable}}\right\}\right..$$

\begin{lem} 
\label{le:rad}
Let $X,Y$ be a pair of representations.
\begin{enumerate}
\item $\Rad(X,Y)$ is a subspace of $\Hom(X,Y)$. 
\item $\Rad(X,Y_1\oplus Y_2)=\Rad(X,Y_1)\oplus\Rad(X,Y_2)$.
\item $\Rad(X_1\oplus X_2,Y)=\Rad(X_1,Y)\oplus\Rad(X_2,Y)$.
\item If $X$ and $Y$ are indecomposable, then
$\Hom(X,Y)\setminus\Rad(X,Y)$ equals the set of isomorphisms $X\to Y$.
\end{enumerate}
\end{lem}
\begin{proof}
(1) Let $\p_1,\p_2\in\Rad(X,Y)$. Choose $\s\in\Hom(Z,X)$ and
$\t\in\Hom(Y,Z)$ with $Z$ indecomposable. Then $\t\p_1\s$ and
$\t\p_2\s$ are non-invertible, and therefore
$\t(\p_1+\p_2)\s=\t\p_1\s+\t\p_2\s$ is non-invertible, since $\End(Z)$
is local by Proposition~\ref{pr:fitting}. Thus $\p_1+\p_2$ belongs to $\Rad(X,Y)$.

(2) Let $Y=Y_1\oplus Y_2$ and $\p=(\p_i)\in\Hom(X,Y)$ with
$\p_i\in\Hom(X,Y_i)$ for $i=1,2$. Choose $\s\in\Hom(Z,X)$ and
$\t=(\t_i)\in\Hom(Y,Z)$ with $Z$ indecomposable and
$\t_i\in\Hom(Y_i,Z)$ for $i=1,2$.  Then
$\t\p\s=\t_1\p_1\s+\t_2\p_2\s$.  

If $\p_i\in \Rad(X,Y_i)$ for $i=1,2$, then $\t_i\p_i\s$ is
non-invertible for $i=1,2$, and therefore $\t\p\s$ is non-invertible,
since $\End(Z)$ is local by Proposition~\ref{pr:fitting}.  Thus $\p$
belongs to $\Rad(X,Y)$. Conversely, let $\p\in\Rad(X,Y)$ and fix
$i\in\{1,2\}$. Then $\p_i\in\Rad(X,Y_i)$ because we can put $\t_j=0$
for $j\neq i$ and have that $\t_i\p_i\s=\t\p\s$ is non-invertible.

(3) Analogous to part (2).

(4) Let $\p\in\Hom(X,Y)\setminus\Rad(X,Y)$. Choose $\s\in\Hom(Z,X)$ and
$\t\in\Hom(Y,Z)$ with $Z$ indecomposable such that $\t\p\s$ is
invertible.  Then $\s$ is invertible because $X$ is indecomposable,
and $\t$ is invertible because $Y$ is indecomposable. Thus $\p$ is
invertible. It is clear that an isomorphism $X\to Y$ does not belong to $\Rad(X,Y)$.
\end{proof}

\subsection{The Krull-Remak-Schmidt theorem}

\begin{thm}
Let $X$ be a finite dimensional representation. Then there exists a
decomposition $X= X_1^{a_1}\oplus\ldots\oplus X_r^{a_r}$ with the
$X_i$ pairwise non-isomorphic indecomposable representations and each
$a_i\geq 1$. If $X= Y_1^{b_1}\oplus\ldots \oplus Y_s^{b_s}$ is another
such decomposition, then $r=s$ and, after reordering, $X_i\cong Y_i$
and $a_i=b_i$ for $i\le i\le r$.
\end{thm}

\begin{proof}
Induction on dimension shows that $X$ decomposes into a finite direct
sum of indecomposable representations. Suppose that
$X=X_1^{a_1}\oplus\ldots\oplus X_r^{a_r}$ is such a direct sum
decomposition with the $X_i$ pairwise non-isomorphic indecomposable
representations and each $a_i\geq 1$.  Let $Y$ be indecomposable and
consider
$$\frac{\dim\Hom(X,Y)-\dim\Rad(X,Y)}{\dim\Hom(Y,Y)-\dim\Rad(Y,Y)}.$$
We see by Lemma~\ref{le:rad} that this number equals $a_i$ if $Y\cong
X_i$ (there is at most one such $i$) and $0$ otherwise. In particular,
this number is independent of the decomposition.
\end{proof}

The Krull-Remak-Schmidt theorem says that the classification of finite
dimensional representations reduces to the classification of
indecomposable representations. There is a similar statement about
morphisms between representations.  Let $X=X_1\oplus\ldots\oplus X_r$
and $Y=Y_1\oplus\ldots\oplus Y_s$ be two representations with their
decompositions into indecomposable representations.  Then we have
$$\Hom(X,Y)=\bigoplus_{i,j}\Hom(X_i,Y_j)$$ by Lemma~\ref{le:sum}. Thus
each morphism $\p\colon X\to Y$ can be written uniquely as a matrix
$\p=(\p_{ij})$ where each entry $\p_{ij}\colon X_i\to Y_j$ is a
morphism between indecomposable representations.

\begin{exm}
Each representation $X$ of the $n$-subspace quiver
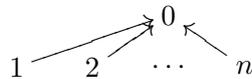
\[\xymatrix@!0 @R=1.8em @C=2.6em {
&&0\\ 1\ar[urr]&2\ar[ur]&\cdots&n\ar[ul]}\] admits a
unique decomposition $X=X'\oplus X(1)\oplus\ldots \oplus X(n)$ such
that $X(i)$ is a direct sum of copies of $S(i)$ for $1\le i\le n$ and
$X'$ is a subspace representation (that is, each linear map $X'_i\to X'_0$ is injective).
\end{exm}

\section{Reflection functors}
In this section we introduce reflection functors. These functors form
our basic tool for classifying representations in terms of their
dimension vectors.

\subsection{Orientations}

A vertex $i$ of $Q$ is called a \emph{sink} (resp.\ \emph{source}) if
there is no arrow in $Q$ starting (resp.\ ending) at $i$.

Given any vertex $i$, the quiver $\s_iQ$ is obtained from $Q$ by
reversing all arrows which start or end at $i$.

An ordering $i_1,\ldots,i_n$ of the vertices of $Q$ is called
\emph{admissible} if for each $p$ the vertex $i_p$ is a sink for
$\s_{i_{p-1}}\ldots\s_{i_1}Q$.  In that case we have
$$\s_{i_n}\ldots\s_{i_1}Q=Q.$$ 

\begin{lem}
There exists an admissible ordering of the vertices of $Q$ if and only
if there are no oriented cycles in $Q$.
\end{lem}
\begin{proof}
We show one implication by induction on the number of vertices. So
suppose $Q$ has no oriented cycles and let $i_n$ be the starting vertex
of a path of maximal length. Then $i_n$ is a source and we remove it
from $Q$. There is an admissible ordering $i_1,\ldots,i_{n-1}$ of the
remaining vertices and we get an admissible ordering
$i_1,\ldots,i_{n}$ of the vertices of $Q$.
\end{proof}

\subsection{The Euler form}
Let $n=\card Q_0$. The \emph{Euler form} is the bilinear form
$$\langle-,-\rangle\colon\bbZ^n\times\bbZ^n\lto\bbZ\quad\text{with}\quad\langle
x,y\rangle=\sum_{i\in Q_0}x_iy_i-\sum_{\a\in Q_1}x_{s(\a)}y_{t(\a)}.$$
We obtain on $\bbZ^n$ a \emph{symmetric bilinear form} by defining
\begin{equation}
\label{eq:bilform}(x,y)=\langle x,y\rangle+\langle y,x\rangle.
\end{equation}

Suppose that $Q$ has no \emph{loops} (that is, arrows from a vertex to
itself).  The \emph{reflection} with respect to a vertex $i$ is by
definition the map 
$$\s_i\colon\bbZ^n\lto\bbZ^n\quad\text{with}\quad
\s_i(x)=x-\frac{2(x,e_i)}{(e_i,e_i)}e_i$$ where $e_i$ is the $i$th
coordinate vector. It is easily checked that the $\s_i$ are
automorphisms of order two preserving the bilinear form $(-,-)$.

For the set $\bbZ^n$ we use the \emph{partial order} which is defined as follows: $$x\le y\quad
\Longleftrightarrow \quad x_i\le y_i\quad\text{for all}\quad i.$$

\subsection{Reflection functors}

Given a pair of quivers $Q$ and $Q'$, a \emph{functor}
$F\colon\Rep(Q,k)\to\Rep(Q',k)$ is an assignment such that
$F\id_X=\id_{FX}$ for every representation $X$ and
$F(\psi\p)=(F\psi)(F\p)$ for every pair $\p\colon X\to Y$ and
$\psi\colon Y\to Z$ of morphisms.

Let $i$ be a vertex of $Q$. We define a pair of \emph{reflection
functors} 
$$S^+_i,S^-_i\colon\Rep(Q,k)\lto\Rep(\s_iQ,k).$$
To this end fix representations $X,X'$ of
$Q$ and a morphism $\p\colon X\to X'$.

(1) If the vertex $i$ is a sink of $Q$, then we construct
$S^+_i$ as follows. We define
$S^+_iX=Y$ by letting $Y_j=X_j$ for a vertex $j\neq i$, and letting
$Y_i$ be the kernel of the map $\xi=(X_\a)$ in the following sequence
$$Y_i\stackrel{\check\xi}\lto\bigoplus_{\substack{\a\in Q_1\\
t(\a)=i}}X_{s(\a)}\stackrel{\xi}\lto X_i$$ where $\check\xi$ denotes
the inclusion map of the kernel.  For an arrow $\a$ in $Q$, let
$Y_\a=X_\a$ if $t(\a)\neq i$, and $Y_\a\colon Y_i\to
X_{s(\a)}=Y_{s(\a)}$ be the map $\check\xi$ followed by the canonical
projection onto $X_{s(\a)}$ if $t(\a)=i$. For the morphism
$S^+_i\p=\psi$, let $\psi_j=\p_j$ if $j\neq i$ and let $\psi_i\colon
Y_i\to Y_i'$ be the restriction of the map
$$(\p_{s(\a)})\colon\bigoplus_{\substack{\a\in Q_1\\
t(\a)=i}}X_{s(\a)}\lto \bigoplus_{\substack{\a\in Q_1\\
t(\a)=i}}X'_{s(\a)}.$$

(2) If the vertex $i$ is a source of $Q$, then we construct dually
$S^-_i$ as follows.  We define $S^-_iX=Y$ by letting $Y_j=X_j$ for a
vertex $j\neq i$, and letting $Y_i$ be the cokernel of the map
$\xi=(X_\a)$ in the following sequence
$$X_i\stackrel{\xi}\lto\bigoplus_{\substack{\a\in Q_1\\
s(\a)=i}}X_{t(\a)}\stackrel{\hat\xi}\lto Y_i$$ where $\hat\xi$
denotes the canonical map onto the cokernel. For an arrow $\a$ in $Q$,
let $Y_\a=X_\a$ if $s(\a)\neq i$, and $Y_\a\colon
Y_{t(\a)}=X_{t(\a)}\to Y_i$ be the restriction of $\hat\xi$ to
$X_{t(\a)}$ if $s(\a)=i$.  For the morphism $S^-_i\p=\psi$, let
$\psi_j=\p_j$ if $j\neq i$ and let $\psi_i\colon Y_i\to Y_i'$ be the
map which is induced by
$$(\p_{t(\a)})\colon\bigoplus_{\substack{\a\in Q_1\\
s(\a)=i}}X_{t(\a)}\lto \bigoplus_{\substack{\a\in Q_1\\
s(\a)=i}}X'_{t(\a)}.$$

(3) Let $i$ be a sink of $Q$. Then we define a
natural monomorphism
\begin{equation}
\label{eq:iota}\iota_iX\colon S^-_iS^+_iX\lto X
\end{equation} 
by letting $(\iota_iX)_j=\id_{X_j}$ for a vertex $j\neq i$, and
letting $(\iota_iX)_i$ be the canonical map
$$(S^-_iS^+_iX)_i=\Coker\check\xi\cong\Im\xi\lto X_i.$$

(4) Let $i$ be a source of $Q$. Then we define a
natural epimorphism
\begin{equation}
\pi_iX\colon X\lto S^+_iS^-_iX
\end{equation} 
by letting $(\pi_iX)_j=\id_{X_j}$ for a vertex $j\neq i$, and letting
$(\pi_iX)_i$ be the canonical map
$$X_i\lto\Im\xi\cong\Ker\hat\xi= (S^+_iS^-_iX)_i.$$

\begin{lem}
\label{le:refl0}
\pushQED{\qed}
For each vertex $i$, $S_i^+$ and $S_i^-$ are functors.\qedhere
\end{lem}

\begin{lem}
\label{le:refl1}
Let $X,X'$ be representations of $Q$ and $i$ be a vertex.
\begin{enumerate}
\item $S_i^\pm(X\oplus X')=S_i^\pm X\oplus S_i^\pm X'$.
\item $X=(S^-_iS^+_iX)\oplus \Coker\iota_iX$ and
$X=(S^+_iS^-_iX)\oplus \Ker\pi_iX$.
\item If $\Coker \iota_iX=0$, then $\dim S^+_iX=\s_i(\dim X)$.
\item If $\Ker \pi_iX=0$, then $\dim S^-_iX=\s_i(\dim X)$.
\end{enumerate}
\end{lem}
\begin{proof} 
(1) Use that $S_i^\pm$ is a functor satisfying
$S^\pm_i(\p+\psi)=S^\pm_i\p+S^\pm_i\psi$ for any pair of parallel
morphisms $\p,\psi$.

(2) The canonical map $\r_i'\colon X_i\to \Coker \xi$ has a section
    $\r_i\colon \Coker\xi\to X_i$, that is,
    $\r_i'\r_i=\id_{\Coker\xi}$. This gives a morphism $\r\colon
    \Coker\iota_iX\to X$ if we put $\r_j=0$ for $j\neq i$. It is
    clear that $\iota_iX\colon S^-_iS^+_iX\to X$ and $\r\colon
    \Coker\iota_iX\to X$ give a direct sum decomposition of $X$.
    The proof for $X=(S^+_iS^-_iX)\oplus \Ker\pi_iX$ is similar.

(3) If $\Coker \iota_iX=0$, then we have
$$\dim Y_i=\sum_{\substack{\a\in Q_1\\ t(\a)=i}}\dim X_{s(\a)}-\dim
X_i$$ and $\dim Y_j=\dim X_j$ for $j\neq i$. Thus $\dim Y=\s_i(\dim
X)$.  The proof of (4) is similar.
\end{proof}

Note that the representations $\Coker\iota_iX$ and $\Ker\pi_iX$
are concentrated at the vertex $i$. Thus they are direct sums of
copies of the simple representation $S(i)$.

\begin{lem}
\label{le:refl2}
Let $i$ be a sink and $X$ an indecomposable representation of
$Q$. Then the following are equivalent:
\begin{enumerate}
\item $X\not\cong S(i)$.
\item $S^+_iX$ is indecomposable.
\item $S^+_iX\neq 0$.
\item $S^-_iS^+_iX\cong X$.
\item The map $(X_\a)\colon \bigoplus_{\substack{\a\in Q_1\\
      t(\a)=i}}X_{s(\a)}\to X_i$ is an epimorphism.
\item $\s_i(\dim X)>0$.
\item $\dim S^+_iX=\s_i(\dim X)$.
\end{enumerate}
\end{lem}
\begin{proof}
Apply  Lemma~\ref{le:refl1}.
\end{proof}

\begin{rem} 
There is an analogue of Lemma~\ref{le:refl2} for a source of $Q$ and
the corresponding functor $S_i^-$.
\end{rem}

The following theorem is a consequence and summarises the basic
properties of the reflection functors.

\begin{thm}
The functors $S^+_i$ and $S^-_i$ induce mutually inverse bijections
between the isomorphism classes of indecomposable representations of
$Q$ and the isomorphism classes of indecomposable representations of
$\s_iQ$, with the exception of the simple representation $S(i)$, which
is annihilated by these functors. Moreover, $\dim S_i^\pm X=\s_i(\dim
X)$ for every indecomposable representation $X$ not isomorphic to
$S(i)$.
\end{thm}

\subsection{Coxeter functors}

Let $Q$ be a quiver without oriented cycles and let $i_1,\ldots,i_n$
be an admissible ordering of the vertices of $Q$.  
The \emph{Coxeter functor} with respect to this ordering is the
functor
$$C^+=S_{i_n}^+\ldots S_{i_1}^+\colon\Rep(Q,k)\lto\Rep(Q,k).$$
We also define
$$C^-=S_{i_1}^-\ldots S_{i_n}^-\colon\Rep(Q,k)\lto\Rep(Q,k).$$
For $r\in\bbZ$, we write
$$C^r=\begin{cases} (C^+)^r&\text{if }r> 0,\\ \Id&\text{if }
r=0,\\(C^-)^{-r}&\text{if }r<0.
\end{cases}$$
Note that in general $C^rC^s\ne C^{r+s}$.

\begin{lem}
The functors $C^+$ and $C^-$ do not depend on the choice of the
admissible ordering of the vertices of $Q$.
\end{lem}
\begin{proof}
First observe that $S_i^+S_j^+=S_j^+S_i^+$ if $i$ and $j$ are sinks
with respect to some orientation and if both are not joined by an
arrow.  Now fix two admissible orderings $i_1,\ldots,i_n$ and
$i'_1,\ldots,i'_n$ of the vertices of $Q$. Let $i_1=i'_m$. Then
$i'_1,\ldots,i'_{m-1}$ are not joined to $i_1$ by an arrow. Therefore
$$S^+_{i'_m}\ldots S^+_{i'_1}=S^+_{i'_{m-1}}\ldots
S^+_{i'_1}S^+_{i_1}.$$ Applying a similar argument for $i_2$, then
$i_3$, and so on, we obtain
\begin{equation*}
S^+_{i'_n}\ldots S^+_{i'_1}=S^+_{i_{n}}\ldots
S^+_{i_1}.\qedhere
\end{equation*}
\end{proof}

For simplicity we assume in the following that $Q_0=\{1,\ldots,n\}$
with $1,\ldots,n$ an admissible ordering.

\begin{lem}
\label{le:dimproj}
Let $i$ be a vertex.
\begin{enumerate}
\item $\dim P(i)=\s_1\ldots\s_{i-1}(e_i)$ and $\dim I(i)=\s_n\ldots\s_{i+1}(e_i)$.
\item $P(i)\cong S^-_1\ldots S^-_{i-1}S(i)$ and $I(i)\cong S^+_n\ldots S^+_{i+1}S(i)$ 
\end{enumerate}
\end{lem}
\begin{proof}
We provide the proof for $P(i)$; the proof for $I(i)$ is dual.
 
(1) For $0\le l < i$, one shows by induction that
\begin{equation}\label{eq:dimP}
\s_{i-l}\ldots\s_{i-1}(e_i)=\sum_{j=0}^l\card Q(i,i-j)e_{i-j}
\end{equation}
in $\bbZ^n$. For $l=i-1$, we then obtain
$\s_{1}\ldots\s_{i-1}(e_i)=\dim P(i)$ because there are no paths from
$i$ to $j$ in $Q$ for $j>i$.

(2) We use the first part, in particular \eqref{eq:dimP}, and apply Lemma~\ref{le:refl2}. An
induction yields for $0\le l< i$ that 
$$\dim S^+_{l}\ldots S^+_1P(i)=\s_{l+1}\ldots\s_{i-1}(e_i).$$ Thus
$S^+_{i-1}\ldots S^+_1P(i)\cong S(i)$ and therefore $P(i)\cong
S^-_1\ldots S^-_{i-1}S(i)$.
\end{proof}

\begin{prop}
\label{pr:proj}
Let $X$ be an indecomposable representation of $Q$.
\begin{enumerate}
\item $C^+X=0$ if $X\cong P(i)$ for some vertex $i$, and $C^-C^+X\cong X$ otherwise.
\item $C^-X=0$ if $X\cong I(i)$ for some vertex $i$, and $C^+C^-X\cong X$ otherwise.
\end{enumerate}
\end{prop}
\begin{proof}
(1) We have $P(i)\cong S^-_1\ldots S^-_{i-1}S(i)$ by
Lemma~\ref{le:dimproj}. Now apply reflection functors and use Lemma~\ref{le:refl2} to obtain
$$C^+P(i)=S^+_n\ldots S^+_iS(i)=0.$$ 
If  $X\not\cong S^-_1\ldots S^-_{i-1}S(i)$ for all $i$, then
\[C^-C^+X=S^-_1\ldots S^-_nS^+_n\ldots S^+_1X\cong
S^-_1\ldots S^-_{n-1}S^+_{n-1}\ldots S^+_1X\cong\ldots\cong X,\]
again by Lemma~\ref{le:refl2}. The proof of (2) is analogous.
\end{proof}

\subsection{Preprojective and preinjective representations}
Let $Q$ be a quiver without oriented cycles. We introduce three
classes of representations.

\begin{defn}
Let $X$ be an indecomposable representation of $Q$.
\begin{enumerate}
\item $X$ is \emph{preprojective} if $X\cong C^rP(i)$ for some vertex $i$ and some $r\le 0$.
\item $X$ is \emph{preinjective} if $X\cong C^rI(i)$ for some vertex $i$ and some $r\ge 0$.
\item $X$ is \emph{regular} if $C^rX\neq 0$ for all $r\in\bbZ$.
\end{enumerate}
\end{defn}

Note that $X$ is preprojective if and only if $C^rX=0$ for some $r>0$,
and $X$ is preinjective if and only $C^rX=0$ for some $r<0$. This is
an immediate consequence of Proposition~\ref{pr:proj}.

\begin{prop}
\label{pr:preproj}
An indecomposable representation is preprojective, preinjective or
regular.  Given indecomposable representations $X,Y$ with $X$
preprojective or preinjective, we have  $X\cong Y$ if and only if $\dim X=\dim Y$. Moreover,
\begin{enumerate}
\item $C^rP(i)\cong C^sP(j)\neq 0$ implies $i=j$ and $r=s$;
\item $C^rI(i)\cong C^sI(j)\neq 0$ implies $i=j$ and $r=s$.
\end{enumerate}
\end{prop}
\begin{proof}
The first assertion is an immediate consequence of
Proposition~\ref{pr:proj}. For the rest we use reflection functors, in
particular Lemma~\ref{le:refl2}. Suppose $\dim X=\dim Y$ and let $X$
be preprojective, say $X\cong C^rP(i)$. We know $\dim P(i)$ from
Lemma~\ref{le:dimproj} and have therefore $$\dim
Y=(\s_n\ldots\s_1)^r\s_1\ldots\s_{i-1}(e_i).$$ Using reflection
functors we obtain $S^+_{i-1}\ldots S^+_1C^{-r}(Y)\cong S(i)$ and
this gives
$$Y\cong C^rS^-_1\ldots S^-_{i-1}S(i)\cong C^rP(i)\cong X.$$ The
proof for preinjective $X$ is analogous.

(1) If $C^rP(i)\cong C^sP(j)\neq 0$ then $P(i)\cong C^{s-r}P(j)$ and
therefore $s-r\le 0$ by Propopsition~\ref{pr:proj}. The same argument
gives $r-s\le 0$ whence $r=s$.  We obtain that $P(i)\cong P(j)$ and
this implies $i=j$ by Lemma~\ref{le:proj}. The proof of (2) is
analogous.
\end{proof}

\section{Dynkin and Euclidean diagrams}
\label{se:graphs}
A finite graph arises from a quiver when one forgets the orientation
of its arrows.  In this section we classify finite graphs using
properties of quadratic forms. For graphs of Dynkin or Euclidean type
we study the corresponding root systems.

\subsection{Finite graphs}

Let $\Ga$ be a finite graph with set of vertices $\{1,\ldots,n\}$.
The finite number of edges joining two vertices $i$ and $j$ is denoted by
$d_{ij}=d_{ji}$.  The graph $\Ga$ induces a \emph{symmetric bilinear form}
$$(-,-)\colon\bbZ^n\times\bbZ^n\lto\bbZ\quad\text{with}\quad(e_i,e_j)=\begin{cases}
-d_{ij}&\text{if }i\neq j,\\ 2-2d_{ii}&\text{if }i=j,
\end{cases}$$
where $e_i$ is the $i$th coordinate vector, and  a
\emph{quadratic form} $$q\colon\bbZ^n\lto\bbZ\quad\text{with}\quad
q(x)=\sum_{i=1}^nx_i^2-\sum_{i\le j}d_{ij}x_ix_j.$$ Note that $\Ga$,
$(-,-)$ and $q$ determine each other, since $q(x)=\frac{1}{2}(x,x)$
and $(x,y)=q(x+y)-q(x)-q(y)$.  The \emph{radical} of the form $q$ is
by definition the set
$$\rad q=\{x\in\bbZ^n\mid (x,-)=0\}.$$ 
A vector $x\in\bbZ^n$ is \emph{sincere} if $x_i\neq 0$ for all $i$.

\begin{rem}
Let $Q$ be a quiver whose underlying graph is $\Ga$. Then the
symmetric bilinear form \eqref{eq:bilform}  for $Q$ coincides with the one
defined for $\Ga$.
\end{rem}

\begin{defn}
Let  $q\colon\bbZ^n\to\bbZ$ be a quadratic form.
\begin{enumerate}
\item $q$ is \emph{positive definite} if $q(x)>0$ for all
non-zero $x\in\bbZ^n$.
\item $q$ is \emph{positive semi-definite} if
$q(x)\ge 0$ for all $x\in\bbZ^n$. 
\end{enumerate}
\end{defn}

\begin{lem}
\label{le:semi-def}
Let $\Ga$ be connected and $y\in\bbZ^n$ be a positive radical
vector. Then $y$ is sincere and $q$ is positive semi-definite. For
$x\in\bbZ^n$, we have
$$q(x)=0 \iff x\in\bbQ y\iff x\in \rad q. $$
\end{lem}
\begin{proof}
The assumption on $y$ yields
\begin{equation}\label{eq:delta}
0=(e_i,y)=(2-2d_{ii})y_i-\sum_{j\neq i}d_{ij}y_j\quad\text{for}\quad1\le i\le n.
\end{equation}
If $y_i=0$ then $\sum_{j\neq i}d_{ij}y_j=0$, and since each term is
non-negative we have $y_j=0$ whenever $i$ and $j$ are joined by an
edge.  It follows that $y=0$ since $\Ga$ is connected. This is a
contradiction and therefore $y$ is sincere. The next calculation shows
that $q$ is positive semi-definite. For $x\in\bbZ^n$, we have
\begin{align*}
q(x)
&= \sum_i(2-2d_{ii})y_i\frac{1}{2y_i}x_i^2 - \sum_{i<j}d_{ij}x_ix_j\\
&= \sum_{i\neq j}d_{ij}\frac{y_j}{2y_i}x_i^2 - \sum_{i<j}d_{ij}x_ix_j\\
&= \sum_{i< j}d_{ij}\frac{y_j}{2y_i}x_i^2 - \sum_{i<j}d_{ij}x_ix_j
+\sum_{i< j}d_{ij}\frac{y_i}{2y_j}x_j^2 \\
&= \sum_{i< j}d_{ij}\frac{y_iy_j}{2}\Big(\frac{x_i}{y_i}-\frac{x_j}{y_j}\Big)^2\ge 0,
\end{align*}
where the second equality follows from \eqref{eq:delta}. If $q(x)=0$
then $\frac{x_i}{y_i}=\frac{x_j}{y_j}$ whenever there is an edge
joining $i$ and $j$. Thus $x\in\bbQ y$ since $\Ga$ is connected. If
$x\in\bbQ y$ then $x\in\rad q$ since $y\in\rad q$ by
assumption. Finally, $x\in\rad q$ implies $q(x)=0$.
\end{proof}

\subsection{The classification}
We list the Dynkin and Euclidean diagrams.
\medskip

\noindent
\emph{Dynkin diagrams} (with $n$ vertices).\\
$
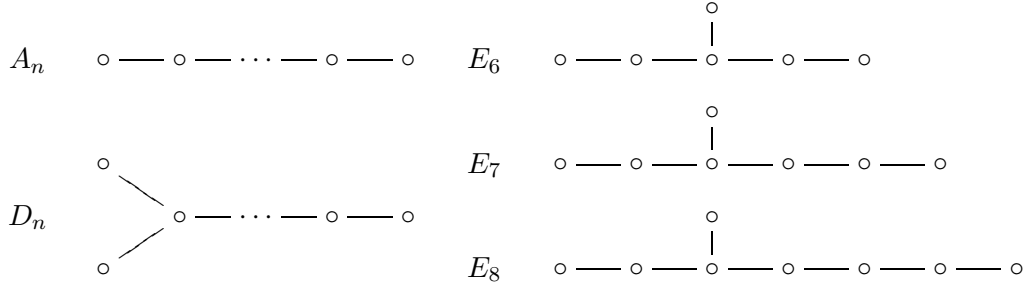
\xymatrix @!0 @R=1.8em @C=2.6em 
{
&&&&&&&&&\circ\ar@{-}[d]\\
A_n&\circ\ar@{-}[r]&\circ\ar@{-}[r]&\cdots\ar@{-}[r]&\circ\ar@{-}[r]&\circ&
E_6&\circ\ar@{-}[r]&\circ\ar@{-}[r]&\circ\ar@{-}[r]&\circ\ar@{-}[r]&\circ\\
&&&&&&&&&\circ\ar@{-}[d]\\
&\circ\ar@{-}[rd]&&&&&
E_7&\circ\ar@{-}[r]&\circ\ar@{-}[r]&\circ\ar@{-}[r]&\circ\ar@{-}[r]&\circ
\ar@{-}[r]&\circ\\
D_n&&\circ\ar@{-}[r]&\cdots\ar@{-}[r]&\circ\ar@{-}[r]&\circ&&&&\circ\ar@{-}[d]\\
&\circ\ar@{-}[ru]&&&&&E_8&\circ\ar@{-}[r]&\circ\ar@{-}[r]&\circ\ar@{-}[r]&\circ\ar@{-}[r]&
\circ\ar@{-}[r]&\circ\ar@{-}[r]&\circ
}$ 

\medskip

\noindent
\emph{Euclidean diagrams} (with $n=m+1$ vertices). Each vertex $i$ is
marked with the value $\d_i$ of a vector $\d\in\bbZ^n$. Let $m\ge 0$
for $\tilde A_m$ and $m\ge 4$ for $\tilde D_m$.\\ 
$
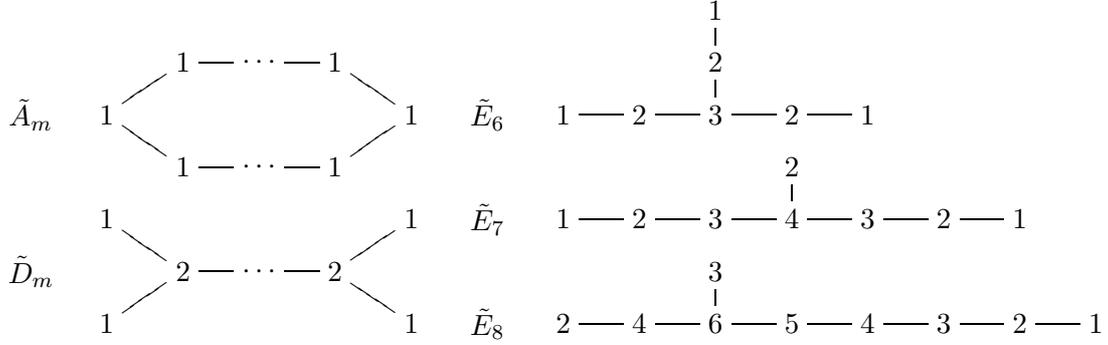
\xymatrix @!0 @R=1.8em @C=2.6em 
{ &&&&&&&&&1\ar@{-}[d]\\
&&1\ar@{-}[ld]\ar@{-}[r]&\cdots\ar@{-}[r]&1\ar@{-}[rd]&&&&&2\ar@{-}[d]\\
\tilde A_m&1&&&&1&\tilde
E_6&1\ar@{-}[r]&2\ar@{-}[r]&3\ar@{-}[r]&2\ar@{-}[r]&1\\
&&1\ar@{-}[lu]\ar@{-}[r]&\cdots\ar@{-}[r]&1\ar@{-}[ru]&&&&&&2\ar@{-}[d]\\
&1\ar@{-}[rd]&&&&1\ar@{-}[ld]&\tilde
E_7&1\ar@{-}[r]&2\ar@{-}[r]&3\ar@{-}[r]&4\ar@{-}[r]&3
\ar@{-}[r]&2\ar@{-}[r]&1\\ \tilde
D_m&&2\ar@{-}[r]&\cdots\ar@{-}[r]&2&&&&&3\ar@{-}[d]\\
&1\ar@{-}[ru]&&&&1\ar@{-}[lu]&\tilde
E_8&2\ar@{-}[r]&4\ar@{-}[r]&6\ar@{-}[r]&5\ar@{-}[r]&
4\ar@{-}[r]&3\ar@{-}[r]&2\ar@{-}[r]&1 }$

\begin{thm} 
Let $\Ga$ be a connected graph and $q$ the corresponding quadratic form.
\begin{enumerate}
\item $\Ga$ is a Dynkin diagram if and only if $q$ is positive
definite.
\item $\Ga$ is a Euclidean diagram if and only if $q$ is positive
semi-definite but not positive definite. In that case there is a unique positive
vector $\d\in\bbZ^n$ with $\rad q=\bbZ \d$. 
\end{enumerate}
\end{thm}
\begin{proof}
We proceed in three steps.

Step 1. \emph{If $\Ga$ is Euclidean then $q$ is positive semi-definite
  and $\rad q=\bbZ \d$.}  The assertion follows from
  Lemma~\ref{le:semi-def} once we have shown that $\d$ is a radical
  vector. This is done by inspection. If $\Ga$ has no loops or
  multiple edges we need to check that
$$0=(e_i,\d)=2\d_i-\sum^n_{\substack{j=1\\ d_{ij}\neq 0}}\d_j\quad\text{for}\quad 1\le i\le n.$$
Finally, since some $\d_i=1$ we have $$\rad q=\bbQ\d\cap\bbZ^n=\bbZ\d.$$

Step 2. \emph{If $\Ga$ is Dynkin then $q$ is positive definite.} This
    follows from the first part because there exists a Euclidean diagram
    $\tilde\Ga$ such that $\Ga$ is obtained by deleting some vertex
    $e$. For $\tilde\Ga$ we have $q(x)>0$ for every non-zero vector
    $x$ with $x_e=0$.

Step 3. \emph{If $\Ga$ is not Dynkin or Euclidean then $q(x)<0$ for some
    $x\in\bbZ^n$.} It is not difficult to find a Euclidean subgraph
    $\Ga'$ with radical vector $\d$. Put $x=\d$ if the vertices of
    $\Ga'$ and $\Ga$ coincide. Otherwise let $i$ be a vertex of
    $\Ga\setminus \Ga'$, but connected with $\Ga'$ by an edge, and take
    $x=2\d+e_i$.

\end{proof}
\subsection{Roots}
Let $\Ga$ be a diagram that is Dynkin or Euclidean.  We define
$$\De=\{x\in\bbZ^n\mid q(x)\le 1\}$$
and a non-zero element of $\De$ is called \emph{root}.

\begin{prop}
\label{pr:root}
Let $\Ga$ be Dynkin or Euclidean.
\begin{enumerate}
\item Each $e_i$ is a root.
\item If $x\in\De$ and $y\in \rad q$, then $-x,x+y\in\De$.
\item Every  root is positive or negative.
\item If $\Ga$ is Euclidean then $\De/{\rad q}$ is finite.
\item If $\Ga$ is Dynkin then $\De$ is finite.
\end{enumerate}
\end{prop}
\begin{proof}
(1) Clear.

(2) We have $q(y\pm x)=q(y)+q(x)\pm (y,x)=q(x)$.

(3) Let $x$ be a root and write $x=x^+-x^-$ where $x^+,x^-\ge
    0$ and both have disjoint support. Then $(x^+,x^-)\le 0$ and
$$1\ge q(x)=q(x^+)+q(x^-)-(x^+,x^-)\ge q(x^+)+q(x^-)\ge 0.$$ This implies
$q(x^+)=0$ or $q(x^-)=0$. Thus one of $x^+$ and $x^-$ is sincere if we
assume that both vectors are non-zero. This is a contradiction and
therefore $x$ is positive or negative.

(4) Fix a vertex $e$. If $x$ is a root with $x_e=0$, then $\d-x$ and
    $\d+x$ are positive at $e$. Both vectors are roots by (2) and
    therefore positive by (3). Thus
$$\{x\in\De\mid x_e=0\}\subseteq\{x\in\bbZ^n\mid-\d\le x\le \d\}$$
which is a finite set. If $x\in\De$ then $x-x_e\d$ belongs to the finite set 
$\{x\in\De\mid x_e=0\}$.

(5) There exists a Euclidean diagram $\tilde\Ga$ such that $\Ga$ is
    obtained by deleting some vertex $e$. A root $x$ of $\Ga$ can be
    viewed as a root for $\tilde\Ga$ with $x_e=0$. Thus the result
    follows from (4).
\end{proof}

\begin{lem}
\label{le:root}
Let $Q$ be a quiver whose underlying graph is Dynkin or Euclidean.  If
$x$ is a positive root and $\s_i(x)$ is not positive, then $x=e_i$.
\end{lem}
\begin{proof}
The root $\s_i(x)$ is not positive by assumption and therefore
negative by Proposition~\ref{pr:root}. For each vertex $j\neq i$, we
have $\s_i(x)_j=x_j$ and therefore $x_j=0$. Thus $x=e_i$.
\end{proof}

\subsection{The Coxeter transformation}
Let $Q$ be a quiver without oriented cycles and fix an admissible
ordering $i_1,\ldots,i_n$ of its vertices. The automorphism
$$c\colon\bbZ^n\lto\bbZ^n\quad\text{with}\quad c(x)=\s_{i_n}\ldots\s_{i_1}(x)$$
is called \emph{Coxeter transformation}. The next lemma shows that
$c$ does not depend on the admissible numbering of the vertices.

\begin{lem}
\label{le:basis}
\begin{enumerate}
\item $c(\dim P(i))=-\dim I(i)$ for every vertex $i$. 
\item $\{\dim P(i)\mid i\in Q_0\}$ and $\{\dim I(i)\mid i\in Q_0\}$
form two bases of $\bbZ^n$.
\end{enumerate}
\end{lem}
\begin{proof}
(1) We simplify the labeling of the vertices and let $i_j=j$ for $1\le
    j\le n$. Then $\dim P(i)=\s_{1}\ldots\s_{i-1}(e_i)$ by
    Lemma~\ref{le:dimproj} and we get $$c(\dim
    P(i))=c\s_{1}\ldots\s_{i-1}(e_i)=\s_n\ldots\s_i(e_i)=-\s_n\ldots\s_{i+1}(e_i)=-\dim
    I(i).$$

(2) We have 
\begin{equation*}e_i=\dim P(i)-\sum_{\substack{\a\in Q_1\\ s(\a)=i}}\dim P(t(\a))
=\dim I(i)-\sum_{\substack{\a\in Q_1\\ t(\a)=i}}\dim I(s(\a)).\qedhere
\end{equation*}
\end{proof}

\begin{lem}
\label{le:gamma}
Let $x,y\in\bbZ^n$.
\begin{enumerate}
\item $\langle \dim P(i),x\rangle=x_i=\langle x,\dim I(i)\rangle$ for each vertex $i$.
\item $\langle x,y\rangle=-\langle y,c(x)\rangle=\langle c(x),c(y)\rangle$.
\end{enumerate}
\end{lem}
\begin{proof}
(1) It is sufficient to check this for each standard basis vector $x=e_j$.

(2) It is sufficient to check this for $x=\dim P(j)$ where $j$ runs
    through all vertices; see Lemma~\ref{le:basis}. Using (1) and
    Lemma~\ref{le:basis} we get
\begin{equation*}\langle\dim P(j),y\rangle=\langle y,\dim I(j)\rangle=\langle
y,-c(\dim P(j))\rangle.\qedhere
\end{equation*}
\end{proof}

\begin{lem}
\label{le:coxtra}
Let $x\in\bbZ^n$. Then $c(x)=x$ if and only if $x\in \rad q$.
\end{lem}
\begin{proof} 
We have $c(x)=x$ iff $x_i=c(x)_i=\s_i(x)_i$ for all $i$ iff
$(x,e_i)=0$ for all $i$ iff $(x,-)=0$.
\end{proof}

From now on we assume that the underlying graph of $Q$ is Dynkin or
Euclidean.  The map $c$ induces a permutation of the finite set
$\De/{\rad q}$. Thus $c^h$ is the identity on $\De/{\rad q}$ for some $h>0$. In
fact, $c^h$ is the identity on $\bbZ^n/{\rad q}$ since $e_i\in\De$ for all
$i$.

\begin{lem}
\label{le:coxtraD}
Let $Q$ be of Dynkin type and $x\in\bbZ^n$.
Then there exists $r\ge 0$
such that $c^r(x)$ is not positive.
\end{lem}
\begin{proof} 
The vector $y=\sum_{r=0}^{h-1}c^r(x)$ is fixed by $c$, and hence
$y=0$ by Lemma~\ref{le:coxtra}. Thus $c^r(x)$ is not positive for
some $r\ge 0$.
\end{proof}

\begin{lem}
\label{le:coxtraE}
Let $Q$ be of Euclidean type and $x\in\bbZ^n$. 
\begin{enumerate}
\item If $c^r(x)>0$ for all $r\in\bbZ$ then $c^h(x)=x$.
\item If $c^h(x)=x$ then $\langle\d,x\rangle=0$.
\end{enumerate}
\end{lem}
\begin{proof} 
(1) Suppose $c^h(x)=x+v$ for some non-zero $v\in \rad q$. An induction
shows that $c^{lh}(x)=x+lv$ for all $l\in\bbZ$. Thus one finds $r$
such that $c^r(x)$ is not positive since the vector $v$ is sincere
and positive or negative.

(2) The vector $y=\sum_{r=0}^{h-1}c^r(x)$ is fixed by $c$, and hence
$y\in\bbZ\d$ by Lemma~\ref{le:coxtra}. Now
$$0=\langle\d,y\rangle=\sum_{r=0}^{h-1}\langle\d,c^r(x)\rangle=h\langle\d,x\rangle$$
since $c$ preserves the Euler form by Lemma~\ref{le:gamma}. Thus
$\langle\d,x\rangle=0$.
\end{proof}

\section{Finite representation type}

In this section we prove Gabriel's theorem. For quivers of Dynkin or
Euclidean type we classify indecomposable representations in terms of
their dimension vectors.

\subsection{The Dynkin case}

\begin{thm}[Gabriel]
\label{th:Gabriel}
Let $Q$ be a quiver whose underlying graph is a Dynkin diagram. Then
the assignment $X\mapsto\dim X$ induces a bijection between the
isomorphism classes of indecomposable representations of $Q$ and the
positive roots corresponding to the diagram of $Q$. In particular,
there are only finitely many isomorphism classes of indecomposable
representations.
\end{thm}
\begin{proof}
Choose an admissible ordering $i_1,\ldots,i_n$ for the vertices of
$Q$. Suppose $X$ is an indecomposable representation of $Q$ with
dimension vector $x=\dim X$. Let 
$$\t=\s_{i_s}\ldots\s_{i_1}(\s_{i_n}\ldots\s_{i_1})^r$$ be the
shortest expression such that $\t(x)$ is not positive, which exists
by Lemma~\ref{le:coxtraD}. Now we apply the reflection functors and
use Lemma~\ref{le:refl2} to obtain
$$S^+_{i_{s-1}}\ldots S^+_{i_1}(S^+_{i_n}\ldots S^+_{i_1})^rX\cong
S(i_{s})$$ and therefore
$$X\cong (S^-_{i_1}\ldots S^-_{i_n})^rS^-_{i_{1}}\ldots
S^-_{i_{s-1}}S(i_{s}).$$ Thus
$$\dim
X=(\s_{i_1}\ldots\s_{i_n})^r\s_{i_1}\ldots\s_{i_{s-1}}(e_{i_{s}})$$ is
a positive root. The same argument shows for another indecomposable
representation $X'$ with $\dim X'=\dim X$ that $X'\cong X$.

Now suppose $x$ is a positive root.  Let
$$\t=\s_{i_s}\ldots\s_{i_1}(\s_{i_n}\ldots\s_{i_1})^r$$ be the
shortest expression such that $\t(x)$ is not positive, which exists by
Lemma~\ref{le:coxtraD}. We infer from Lemma~\ref{le:root} that
$$\s_{i_{s-1}}\ldots\s_{i_1}(\s_{i_n}\ldots\s_{i_1})^r(x)=e_{i_{s}}.$$ Let
$$X=(S^-_{i_1}\ldots S^-_{i_n})^rS^-_{i_{1}}\ldots
S^-_{i_{s-1}}S(i_{s}).$$ Another iterated application of reflection
functors shows that $X$ is indecomposable with
$$\dim X=(\s_{i_1}\ldots \s_{i_n})^r\s_{i_{1}}\ldots
\s_{i_{s-1}}(e_{i_{s}})=x.$$

There are only finitely many isomorphism classes of indecomposable
representations because the number of roots is finite by Proposion~\ref{pr:root}.
\end{proof}

\begin{rem}
The proof of Theorem~\ref{th:Gabriel} shows that for a quiver of
Dynkin type all indecomposable representations are preprojective and
preinjective.  In fact, the proof simplifies a bit if one uses Coxeter
functors.
\end{rem}

\begin{rem}
Given a graph $\Ga$ without loops, the corresponding \emph{Weyl group}
$W(\Ga)$ is the group of automorphisms of $\bbZ^n$ which is generated
by the reflections $\s_i$. If $\Ga$ is a Dynkin diagram then the roots
are precisely the vectors of the form $\t(e_i)$ with $\t\in W(\Ga)$
and $1\le i\le n$. Clearly, each $\t(e_i)$ is a root since the $\s_i$
preserve the quadratic form $q$. Conversely, let $q(x)=1$. Then the
argument used in the proof of Theorem~\ref{th:Gabriel} shows that $x$
is of the form $\t(e_i)$.
\end{rem}

\subsection{The defect}

Let $Q$ be a quiver of Euclidean type. The \emph{defect} of a vector
$x\in\bbZ^n$ is
$$\partial x=\langle\d,x\rangle=-\langle x,\d\rangle.$$ The
\emph{defect} of a representation $X$ is $\partial X=\partial\dim X$.

\begin{prop}
\label{pr:defect}
Let $X$ be an indecomposable representation.
\begin{enumerate}
\item $X$ is preprojective if and only if $\partial X<0$.
\item $X$ is preinjective if and only if $\partial X>0$.
\item $X$ is regular if and only if $\partial X=0$.
\end{enumerate}
\end{prop}
\begin{proof}
First observe that for every representation $X$ with $C^rX\neq 0$ we have
$$\dim C^r X=c^r(\dim X)$$ by Lemma~\ref{le:refl2}. Now let
$X=C^rP(i)$ be preprojective. Then
$$\partial X=-\langle c^r(\dim P(i)),\d\rangle=-\langle\dim
P(i),\d\rangle=-\d_i<0$$ by Lemmas~\ref{le:gamma} and
\ref{le:coxtra}. Similarly preinjectives have positive defect. If $X$
is regular then $\partial X=0$ by Lemma~\ref{le:coxtraE}.
\end{proof}

\subsection{The Euclidean case}

\begin{thm}
\label{th:tame}
Let $Q$ be a quiver without oriented cycles and suppose the underlying
graph is a Euclidean diagram with $n$ vertices. Then the assignment
$X\mapsto \dim X$ induces a bijection between the isomorphism classes
of indecomposable preprojective or preinjective representations of $Q$
and the positive roots with non-zero defect corresponding to the
diagram of $Q$. The preprojective and preinjective indecomposables
form $2n$ countably infinite series $C^{-r}P(i)$ and $C^rI(i)$
($r\in\bbN_0,i\in Q_0$) of pairwise non-isomorphic representations.
\end{thm}
\begin{proof}
Choose an admissible ordering $i_1,\ldots,i_n$ for the vertices of
$Q$. Suppose $X$ is an indecomposable representation of $Q$ with
dimension vector $x=\dim X$. Let $X\cong C^rP(i_s)$ be preprojective.
We have $\dim P(i_s)=\s_{i_1}\ldots\s_{i_{s-1}}(e_{i_s})$ by
Lemma~\ref{le:dimproj} and therefore
$x=c^r\s_{i_1}\ldots\s_{i_{s-1}}(e_{i_s})$ by
Lemma~\ref{le:refl2}. This is a positive root, and $\partial x<0$ by
Proposition~\ref{pr:defect}.  A similar argument works if $X$ is
preinjective. The map $X\mapsto\dim X$ is injective by
Proposition~\ref{pr:preproj}.

Now suppose $x$ is a positive root with $\partial x\neq 0$.  We know
from Lemma~\ref{le:coxtraE} that $c^t(x)$ is not positive for some
$t\in\bbZ$.  Suppose first $t>0$. Then there are $1\le s\le n$ and
$r\ge 0$ such that 
$$\t=\s_{i_s}\ldots\s_{i_1}(\s_{i_n}\ldots\s_{i_1})^{r}$$ is the
shortest expression with $\t(x)$  not positive.
We infer from Lemma~\ref{le:root} that
$$\s_{i_{s-1}}\ldots\s_{i_1}(\s_{i_n}\ldots\s_{i_1})^r(x)=e_{i_{s}}.$$ Let
$$X=(S^-_{i_1}\ldots S^-_{i_n})^rS^-_{i_{1}}\ldots
S^-_{i_{s-1}}S(i_{s}).$$ An iterated application of reflection
functors and Lemma~\ref{le:refl2} shows that $X$ is indecomposable
with
$$\dim X=(\s_{i_1}\ldots \s_{i_n})^r\s_{i_{1}}\ldots
\s_{i_{s-1}}(e_{i_{s}})=x.$$ Moreover, $X$ is preprojective since
$X\cong C^{-r}P(i_s)$ by Lemma~\ref{le:dimproj}. In case $t<0$ a
similar argument gives a preinjective representation $X$ with $\dim X=x$.

Finally we show that the preprojectives and preinjectives form $2n$
countably infinite series of pairwise non-isomorphic
representations. This follows essentially from
Proposition~\ref{pr:preproj}. It remains to show that $C^{-r}P(i)\neq
0$ and $C^rI(i)\neq 0$ for all $r\ge 0$. But this follows from
Poposition~\ref{pr:proj} because $C^{-r}P(i)= 0$ would imply that
$P(i)$ is preinjective, contradicting the fact that preprojectives and
preinjectives have different defect by Proposition~\ref{pr:defect}. A
similar argument works for $C^rI(i)$.
\end{proof}

\begin{prop}
\label{pr:cycle}
Let $Q$ be a quiver of Euclidean type $\tilde A_n$ with $n\ge 0$.
Then there are infinitely many isomorphism classes of indecomposable
representations.
\end{prop}
\begin{proof}
We allow any orientation, in particular an oriented cycle, and fix an
arrow $\a_0$. Now we define for each $p\ge 1$ a representation $X=X(p)$ as
follows. Let $X_i=k^p$ for each vertex $i$, let $X_{\a_0}=J(p,0)$ be the
Jordan block of size $p$ with eigenvalue $0$ and let
$X_\a=\id_{k^p}$ for every arrow $\a\neq\a_0$. Then $\End(X(p))\cong
k[t]/(t^p)$ and therefore $X(p)$ is indecomposable.
\end{proof}

\begin{cor}[Gabriel]
\label{co:gabriel}
Let $Q$ be a connected quiver. Then there are only finitely many
isomorphism classes of indecomposable representations if and only if
the underlying graph is a Dynkin diagram.
\end{cor}
\begin{proof} 
If $Q$ is of Dynkin type then the classification of the indecomposable
representations in Theorem~\ref{th:Gabriel} shows that there are only
finitely many. If $Q$ is not of Dynkin type then $Q$ has a Euclidean
subquiver $Q'$ which has infinitely many indecomposable
representations by Theorem~\ref{th:tame} and
Proposition~\ref{pr:cycle}. Each representation $X$ of $Q'$ can be
extended to a representation of $Q$ by letting $X_i=0$ and $X_\a=0$
for all $i\in Q_0$ and $\a\in Q_1$ not in $Q'$. Thus $Q$ has
infinitely many pairwise non-isomorphic indecomposable
representations.
\end{proof}

\section{Irreducible morphisms}

In this section we investigate the morphisms between two
representations using the concept of an irreducible morphism.  We
think of irreducible morphisms as generators, and in some cases, all
morphisms are sums of compositions of such irreducible morphisms.

\subsection{The radical}
Let $X,Y$ be a pair of representations. Recall from
Section~\ref{se:radical} that the radical $\Rad(X,Y)$ is a subspace of
$\Hom(X,Y)$.  We extend this definition recursively for each $n\ge 0$
as follows. Let $\Rad^0(X,Y)=\Hom(X,Y)$ and for $n> 0$ let
$\Rad^n(X,Y)$ be the set of morphisms $\p\in\Hom(X,Y)$ which admit a
factorisation $\p=\p''\p'$ with $\p'\in\Rad(X,Z)$ and
$\p''\in\Rad^{n-1}(Z,Y)$ for some representation $Z$.  
Note that
$\Rad^1(X,Y)=\Rad(X,Y)$.

\begin{lem}\label{le:radn}
Let $X,Y,Z$ be representations and $m,n\ge 0$.
\begin{enumerate}
\item $\Rad^{n+1}(X,Y)$ is a subspace of $\Rad^n(X,Y)$.
\item For each finite set of representations $X_i$ and $Y_j$, we have
\begin{equation*}
\label{eq:rad}
\Rad^n(\bigoplus_{i}X_i,\bigoplus_{j}Y_j)=\bigoplus_{i,j}\Rad^n(X_i,Y_j).
\end{equation*}
\item If $\p\in\Rad^n(X,Y)$ and $\psi\in\Rad^{m}(Y,Z)$, then
$\psi\p\in\Rad^{n+m}(X,Z)$.
\end{enumerate}
\end{lem}
\begin{proof}
Use Lemma~\ref{le:rad}.
\end{proof}

\subsection{Irreducible morphisms}
Fix a morphism $\p\colon X\to Y$ between two representations. 

The morphism $\p$ is called \emph{split monomorphism} if there exists
$\p'\colon Y\to X$ with $\p'\p=\id_X$, and $\p$ is called \emph{split
epimorphism} if there exists $\p''\colon Y\to X$ with $\p\p''=\id_Y$.

The morphism $\p$ is called \emph{irreducible} if $\p$ is neither a
split monomorphism nor a split epimorphism and if for any
factorisation $\p=\p''\p'$ the morphism $\p'$ is a split monomorphism
or the morphism $\p''$ is a split epimorphism.

\begin{lem} 
\label{le:irr_mon_epi}
An irreducible morphism is a monomorphism or an epimorphism.
\end{lem} 
\begin{proof}
Let $\p\colon X\to Y$ be irreducible. Consider the canonical
factorisation $X\to\Im\p\to Y$. If $X\to\Im\p$ is a split monomorphism
then $\p$ is a monomorphism. If $\Im\p\to Y$ is a split epimorphism then
$\p$ is an epimorphism.
\end{proof}

\begin{lem}
\label{le:irr}
Let $\p\colon X\to Y$ be a morphism between two representations.
\begin{enumerate}
\item If $X$ is indecomposable, then $\p\in\Rad(X,Y)$ if and only if
$\p$ is not a split mono.
\item If $Y$ is indecomposable, then $\p\in\Rad(X,Y)$ if and only if
$\p$ is not a split epi.
\item If $X$ and $Y$ are indecomposable, then
$\p\in\Rad^1(X,Y)\setminus\Rad^2(X,Y)$ if and only if $\p$ is
irreducible.
\end{enumerate}
\end{lem}
\begin{proof}
(1) We apply Lemma~\ref{le:rad}. Let $Y=\bigoplus_{i=1}^rY_i$ be a
    decomposition into indecomposable represenations. Then $\Rad(X,Y)=
    \bigoplus_{i=1}^r\Rad(X,Y_i)$ and we have $\p\in\Rad(X,Y)$ if and
    only if each component $\p_i$ belongs to $\Rad(X,Y_i)$. Thus
    $\p\not\in\Rad(X,Y)$ if and only if $\p_i$ is an
    isomorphism for some $i$ if and only if $\p$ is a split monomorphism.

(2) This is the dual statement of (1).

(3) Combine (1) and (2).
\end{proof}
We define for indecomposable representations $X,Y$
$$\Irr(X,Y)=\Rad^1(X,Y)/\Rad^2(X,Y).$$

\begin{rem}
Let $\p\colon X\to Y$ be an irreducible morphism between
indecomposable representations and let $Y'=Y\oplus Y$.  Then the
morphism $(\p,\p)\colon X\to Y'$ is not irreducible but belongs
to $\Rad^1(X,Y')\setminus\Rad^2(X,Y')$.
\end{rem}

\begin{prop}
\label{pr:irr}
Let $X,Y$ be indecomposable representations and $\Rad^n(X,Y)=0$ for
some $n$. Then every non-isomorphism $X\to Y$ is a sum of compositions
of irreducible morphisms between indecomposable representations.
\end{prop}
\begin{proof}
Suppose $\p$ is not irreducible. Because $\p$ is not an isomorphism,
there exists a factorisation $\p=\p''\p'$ with $\p'\in\Rad(X,Z)$ and
$\p''\in\Rad(Z,Y)$ for some representation $Z$ by
Lemma~\ref{le:irr}. Let $Z=\bigoplus_{i=1}^rZ_i$ be a decomposition
into indecomposable representations.  Then
$\p=\sum_{i=1}^r\p_i''\p_i'$ with $\p_i'\in\Rad(X,Z_i)$ and
$\p_i''\in\Rad(Z_i,Y)$ for all $i$. Next we factorise each
non-irreducible $\p_i'$ and each non-irreducible $\p_i''$ into a
composition of two radical morphisms. We continue and obtain in each
step a finite sum $\p=\sum_{i}\p_{in_i}\ldots\p_{i1}$ of compositions
of radical morphisms between indecomposable representations. This
process stops and all $\p_{ij}$ are irreducible because a composition
of $n$ radical morphisms is zero by our assumption on $X$ and $Y$.
\end{proof}

\subsection{The Harada-Sai lemma}

The \emph{length} $\ell(X)$ of a representation $X$ is the maximal number $n$
such that there exists a chain of subrepresentations 
\[0=X_n\subsetneq \ldots\subsetneq X_1\subsetneq X_0=X.\] 
Note that we have $\ell(X)=\sum_{i\in Q_0}\dim X_i$ and $\ell (X)=\ell
(U)+\ell(X/U)$ for every subrepresentation $U\subseteq X$.

\begin{lem}[Harada-Sai]
\label{le:HS}
Let $\p_i\colon X_i\to X_{i+1}$ with $1\le i\le 2^n-1$ be a family of
non-isomorphisms between indecomposable representations satisfying
$\ell(X_i)\le n$ for all $i$. Then $\p_{2^n-1}\ldots\p_1=0$.
\end{lem}
\begin{proof}
We show by induction that the length of the image of
$\p_{2^m-1}\ldots\p_1$ is at most $n-m$. This is clear for $m=1$ since
$\p_1$ is not an isomorphism. Let $\p'=\p_{2^{m-1}-1}\ldots\p_1$,
$\p=\p_{2^{m-1}}$, and $\p''=\p_{2^m-1}\ldots\p_{2^{m-1}+1}$.  By the
inductive hypothesis, the length of $\Im\p'$ and $\Im\p''$ is at most
$n-m+1$. The assertion follows if either is strictly less. So suppose
the images of $\p'$,$\p''$ and $\p''\p\p'$ each has length
$n-m+1$. Then $\Ker (\p''\p)\cap\Im\p'=0$ and
$\Ker\p''\cap\Im(\p\p')=0$. On the other hand, we have
\begin{gather*}\ell (X_{2^{m-1}})=\ell (\Ker (\p''\p))+\ell (\Im(\p''\p))=\ell (\Ker
(\p''\p))+\ell (\Im\p'),\\ \ell (X_{2^{m-1}+1})=\ell (\Ker \p'')+\ell
(\Im\p'')=\ell (\Ker \p'')+\ell (\Im(\p\p')).
\end{gather*}
Therefore $X_{2^{m-1}}=\Ker
(\p''\p)\oplus\Im\p'$ and
$X_{2^{m-1}+1}=\Ker\p''\oplus\Im(\p\p')$. Since each is
indecomposable, $\p''\p$ is a monomorphism and $\p\p'$ is an
epimorphism. Thus $\p$ is an isomorphism, contrary to hypothesis.
\end{proof}

\begin{prop}
\label{pr:rad}
Suppose that the length of every indecomposable
representation is bounded by $n$. Then $\Rad^{2^n-1}(X,Y)=0$ for every pair $X,Y$ of
representations.
\end{prop}
\begin{proof}
We may assume that $X$ and $Y$ are indecomposable, by
Lemma~\ref{le:radn}. An induction shows that any morphism in
$\Rad^r(X,Y)$ can be written as a finite sum
$\sum_{i}\p_{ir}\ldots\p_{i1}$ of compositions of radical morphisms
between indecomposable representations.  Thus $\Rad^{2^n-1}(X,Y)=0$ by
Lemma~\ref{le:HS}.
\end{proof}

\section{Morphisms between preprojective representations}

We give a combinatorial description of all morphisms between
preprojective representations.

Throughout this section $Q$ denotes a quiver without oriented cycles.

\subsection{Irreducible morphisms between indecomposable projectives}
Given a pair $i,j$ of vertices, each arrow $\a\colon i\to j$ induces a
morphism
$$\a^*\colon P(j)\lto P(i)\quad\text{with}\quad \a^*_l=k[Q(\a,l)]\quad
(l\in Q_0).$$ Note that $\a^*$ is a monomorphism since each map
$Q(\a,l)$ is injective. Taking all arrows starting at $i$ (resp.\
ending at $j$) we obtain two morphisms
$$\s(i)\colon\bigoplus_{\a\colon i\to i'}P(i')\xto{(\a^*)}
P(i)\quad\text{and}\quad \t(j)\colon
P(j)\xto{(\a^*)}\bigoplus_{\a\colon j'\to j}P(j').
$$

\begin{lem}
\label{le:max}
The morphism $\s(i)$ is a monomorphism and its image is the unique
maximal subrepresentation of $P(i)$.
\end{lem}
\begin{proof}
This is a direct consequence of the definition of $P(i)$. Note that a
subrepresentation $X\subseteq P(i)$ equals $P(i)$ if and only if the
trivial path $\e_i$ belongs to $X_i$.
\end{proof}

\begin{lem}
\label{le:irr_proj}
For a  morphism $\p\colon X\to P(i)$, the following are equivalent:
\begin{enumerate}
\item $\p\in\Rad(X,P(i))$.
\item $\p$ is not an epimorphism.
\item $\p$ admits a factorisation $\p=\s(i)\p'$.
\end{enumerate}
\end{lem}
\begin{proof}
(1) $\Leftrightarrow$ (2): We have $\p\in\Rad(X,P(i))$ if and only if
$\p$ is not a split epimorphism, by Lemma~\ref{le:irr}. Thus we need
to show that every epimorphism is a split epimorphism. So let $\p$ be
an epimorphism.  Then the map $\p_i$ is surjective and we find $x\in X_i$
with $\p_i(x)=\e_i$. It follows from Lemma~\ref{le:hom} that there is
a morphism $\p'\colon P(i)\to X$ corresponding to $x$. Then we have
$\p\p'=\id_{P(i)}$ and therefore $\p$ is a split epimorphism.

(2) $\Leftrightarrow$ (3): Apply Lemma~\ref{le:max}.
\end{proof}

\begin{lem}
\label{le:irr_proj1}
Let $i$ be a vertex and $X$ an indecomposable representation. 
\begin{enumerate}
\item If $X\to P(i)$ is an irreducible morphism then there is an arrow
$i\to j$ such that $X\cong P(j)$.
\item Suppose $P(i)$ is simple. If $P(i)\to X$ is an irreducible
morphism then there is an arrow $j\to i$ such that $X\cong P(j)$.
\end{enumerate}
\end{lem}
\begin{proof}
(1) Let $\p\colon X\to P(i)$ be irreducible. Then $\p\in\Rad(X,P(i))$
    by Lemma~\ref{le:irr} and we obtain a factorisation $\p=\s(i)\p'$
    by Lemma~\ref{le:irr_proj}.  The morphism $\s(i)$ is not a split
    epimorphism and therefore $\p'$ is a split monomorphism. Thus
    $X\cong P(j)$ for some arrow $i\to j$.

(2) Let $\p\colon P(i)\to X$ be irreducible. We claim that $\p$
factors through $\t(i)$. First observe that $X\not\cong P(i)=S(i)$. This implies
$$\xi\colon\bigoplus_{\a\colon i'\to i}X_{i'}\xto{(X_\a)} X_i$$ is an epimorphism, by
Lemma~\ref{le:refl2}. We obtain the following commutative diagram
$$\xymatrix{\bigoplus_{\a\colon i'\to
i}\Hom(P(i'),X)\ar@{=}[r]\ar[d]^\wr& \Hom(\bigoplus_{\a\colon i'\to
i}P(i'),X)\ar[rr]^-{\Hom(\t(i),X)}&&
\Hom(P(i),X)\ar[d]^\wr\\ \bigoplus_{\a\colon i'\to
i}X_{i'}\ar[rrr]^-\xi&&& X_i }$$ where the vertical maps are taken from
Lemma~\ref{le:hom}. Thus $\Hom(\t(i),X)$ is surjective and therefore
$\p$ factors through $\t(i)$. Let $\p=\p'\t(i)$ be a factorisation.
The morphism $\t(i)$ is not a split monomorphism and therefore $\p'$
is a split epimorphism. Thus $X\cong P(j)$ for some arrow $j\to i$.
\end{proof}

We denote by $Q_1(i,j)$ the set of arrows $i\to j$ in $Q$.
\begin{lem}
\label{le:irr_proj2}
The map sending an arrow $\a\colon i\to j$ to $\a^*$ induces an
isomorphism $$f\colon k[Q_1(i,j)]\xto{\sim}\Irr(P(j),P(i)).$$
\end{lem}
\begin{proof}
The map $f$ is well-defined because for each arrow $\a\colon i\to j$
the morphism $\a^*$ belongs to $\Rad(P(j),P(i))$ and induces therefore
an element of $\Irr(P(j),P(i))$.  

To show that $f$ is an epimorphism choose a radical morphism $\p\colon
P(j)\to P(i)$.  Then $\p$ admits a factorisation $\p=\s(i)\p'$ by
Lemma~\ref{le:irr_proj}.  Let $\p'=(\p'_\a)$ with
$\p'_\a\in\Hom(P(j),P(t(\a)))$.  If $t(\a)=j$ then we have
$\p'_\a=\la_\a\id_{P(j)}$ for some $\la_\a\in k$ because
$\End(P(j))\cong k$ by Lemma~\ref{le:endproj}.  If $t(\a)\neq j$ then
$\p'_\a\in\Rad(P(j),P(t(\a)))$.  It follows that $f(\sum_{\a\in
Q_1(i,j)}\la_\a\a)$ and $\p$ represent the same element in
$\Irr(P(j),P(i))$ since $\p-\sum_{\a\in
Q_1(i,j)}\la_\a\a^*\in\Rad^2(P(j),P(i))$. 

To show that $f$ is a monomorphism choose $x=\sum_{\a\in
Q_1(i,j)}\la_\a\a$ with $f(x)=0$. Thus $\theta=\sum_{\a\in
Q_1(i,j)}\la_\a\a^*$ belongs to $\Rad^2(P(j),P(i))$ and there is a
factorization $\theta=\p\psi$ such that $\p$ and $\psi$ are radical
morphisms.  We obtain a second factorisation $\p=\s(i)\p'$ by
Lemma~\ref{le:irr_proj} and consider $\p'\psi\colon P(j)\to
\bigoplus_{\a\colon i\to i'} P(i')$. Then
$(\p'\psi)_\a=\la_\a\id_{P(j)}$ if $t(\a)=j$, and $(\p'\psi)_\a=0$ if
$t(\a)\neq j$, since 
$$\sum_{\a\in Q_1(i,j)}\la_\a\a^*=\s(i)\p'\psi=\sum_{\a\colon i\to
i'}\a^*(\p'\psi)_\a$$ and $\s(i)$ is a monomorphism. This implies
$\la_\a=0$ for all $\a\colon i\to j$ since $(\p'\psi)_\a$ is a radical
morphism. Thus $f$ is a monomorphism.
\end{proof}

Lemma~\ref{le:irr_proj2} shows that for each arrow $\a$
the morphism $\a^*$ is irreducible.

\subsection{More irreducible morphisms}
We fix an arrow $\a\colon i\to j$ in $Q$ and construct another irreducible morphism 
$$\a_*\colon P(i)\lto C^-P(j)$$ as follows. We may assume that
$1,\ldots,n$ is an admissible numbering of the vertices of $Q$. Let
$\tilde\a\colon j\to i$ be the arrow in $\tilde Q=\s_{i-1}\ldots\s_1
Q$ corresponding to $\a$. This induces the morphism
$\tilde\a^*\colon\tilde P(i)\to\tilde P(j)$ between representations of
$\tilde Q$ and we define $\a_*=S^-_1\ldots S^-_{i-1}\tilde\a^*$. Note
that we can identify $P(i)=S^-_1\ldots S^-_{i-1}\tilde P(i)$ and
$C^-P(j)=S^-_1\ldots S^-_{i-1}\tilde P(j)$ by Lemma~\ref{le:dimproj},
since $\tilde P(i)=S(i)$ and $\tilde P(j)= S^-_i\ldots S^-_{j-1}S(j)$.

\subsection{Reflection functors and morphisms}

\begin{lem} 
\label{le:reflmor}
Let $i$ be a sink and $X,Y$ indecomposable representations not
isomorphic to $S(i)$. Then $S_i^+$ induces  isomorphisms
$$\Rad^n(X,Y)\xto{\sim}\Rad^n(S^+_iX,S^+_iY)\quad\text{for}\quad
n\ge 0.$$ In particular, $S_i^+$ induces  isomorphisms
\[\Hom(X,Y)\xto{\sim}\Hom(S^+_iX,S^+_iY)\quad\text{and}\quad
\Irr(X,Y)\xto{\sim}\Irr(S^+_iX,S^+_iY).\]
\end{lem}
\begin{proof}
We use the natural morphism $\iota_iZ\colon S_i^-S_i^+Z\to Z$
\eqref{eq:iota} which is defined for any representation $Z$; it is a
split monomorphism by Lemma~\ref{le:refl1}.  Thus we can identify
$S^-_iS^+_iX=X$ and $S^-_iS^+_iY=Y$.  Using this identification the
inverse for $\Hom(X,Y)\to\Hom(S^+_iX,S^+_iY)$ sends $\psi\in
\Hom(S^+_iX,S^+_iY)$ to $S^-_i\psi$. Now fix $\p\in\Hom(X,Y)$.
Clearly, $\p$ is an isomorphism if and only if $S^+_i\p$ is an
isomorphism. Thus $S^+_i$ induces a bijection
$$\Rad^1(X,Y)\xto{\sim}\Rad^1(S^+_iX,S^+_iY).$$ Next we suppose
$\p\in\Rad^n(X,Y)$ and $n>1$. Then $\p$ admits a factorisation
$\p=\p''\p'$ with $\p'\in\Rad^1(X,Z)$ and $\p''\in\Rad^{n-1}(Z,Y)$ for
some representation $Z$. We know by induction that
$S^+_i\p'\in\Rad^1(S^+_iX,S^+_iZ)$ and
$\p''\in\Rad^{n-1}(S^+_iZ,S^+_iY)$. Thus
$S_i^+\p\in\Rad^n(S_i^+X,S_i^+Y)$. The same argument shows that
$S^-_i$ maps $\Rad^n(S^+_iX,S^+_iY)$ to $\Rad^n(X,Y)$.  This
establishes for all $n>1$ the isomorphism
\begin{equation*}\Rad^n(X,Y)\xto{\sim}\Rad^n(S^+_iX,S^+_iY).\qedhere
\end{equation*}
\end{proof}

There are a number of consequences.

\begin{prop}\label{pr:coxmor}
Let $X,Y$ be representations having no indecomposable projective
direct summand. Then $C^+$ induces an isomorphism
$\Hom(X,Y)\xto{\sim}\Hom(C^+X,C^+Y)$.
\end{prop}
\begin{proof}
Combine Lemma~\ref{le:reflmor} and Proposition~\ref{pr:proj}.
\end{proof}

\begin{prop}
\label{pr:end}
If $X$ is an indecomposable preprojective or preinjective
representation then $\End(X)\cong k$.
\end{prop}
\begin{proof}
Combine Lemmas~\ref{le:endproj} and \ref{le:reflmor}.
\end{proof}

\begin{prop}
\label{pr:irrpp}
Let $X=C^rP(i)$ and $Y=C^sP(j)$ be two indecomposable preprojective
representations. Then we have
$$\Irr(X,Y)\cong\begin{cases} k[Q_1(j,i)]&\text{if }r=s,\\
k[Q_1(i,j)]&\text{if }r=s+1,\\
0&\text{otherwise}.
\end{cases}$$
The isomorphism sends an arrow $\a\colon j\to i$ to $C^r\a^*$ and
$\b\colon i\to j$ to $C^r\b_*$.
\end{prop}
\begin{proof}
Suppose there exists an irreducible morphisms $X\to Y$. 

Let $r\le s$. Then $C^{r-s}P(i)\neq 0$ and we have
$$\Irr(X,Y)\cong\Irr(C^{r-s}P(i),P(j))\cong k[Q_1(j,i)],$$ where the
first isomorphism follows from Lemma~\ref{le:reflmor}, and the second
follows from Lemmas~\ref{le:irr_proj1} and \ref{le:irr_proj2}. In
particular, $r=s$.

A similar argument works if $r>s$.  Then $C^{s-r}P(j)$ is an
indecomposable and non-projective representation. We may assume that
$1,\ldots,n$ is an admissible numbering of the vertices of $Q$ and
we have
\begin{align*}
\Irr(X,Y)&\cong\Irr(P(i),C^{s-r}P(j)\\
&\cong\Irr(S_{1}^-\ldots S_{i-1}^-S(i),S_1^-\ldots S_{i-1}^-S_i^-\ldots S_n^-
C^{s-r+1}S_1^-\ldots S_{j-1}^-S(j))\\
&\cong\Irr(\tilde P(i),\tilde P(j))\cong k[Q_1(i,j)],
\end{align*}
where $\tilde P(i)=S(i)$ and $\tilde P(j)= S^-_i\ldots
S^-_{j-1}S(j)$ denote the indecomposable projective representations of
$\tilde Q= \s_{i-1}\ldots\s_1Q$ corresponding to $i$ and $j$
respectively. In particular, $r=s+1$.

It is clear that this argument can be reversed. Thus we have a
necessary and sufficient criterion for the existence of irreducible
morphisms $X\to Y$.
\end{proof}

\subsection{Reflection functors and exact sequences}

\begin{lem}\label{le:monoepi}
Given  a vertex $i$ of $Q$, any exact sequence $0\to X'\to X\to
X''\to 0$ induces exact sequences
\begin{gather*}
0\to S_i^+X'\to S_i^+X\to S_i^+X''\to X^+\to 0\quad(\text{$i$ a sink})\\
0\to X^-\to  S_i^-X'\to S_i^-X\to S_i^-X''\to 0\quad(\text{$i$ a source})
\end{gather*}
such that $X^+$ and $X^-$ are direct sums of copies of $S(i)$.
\end{lem}
\begin{proof}
Apply Lemma~\ref{le:snake}.
\end{proof}

\begin{lem}
Suppose $Q$ has no oriented cycles and fix a vertex $i$. Then we have
the following exact sequences:
\begin{gather}
\label{eq:1}0\lto\bigoplus_{\a\colon i\to j}
P(j)\stackrel{(\a^*)}\lto P(i)\lto  S(i)\lto 0\\
\label{eq:2}0\lto P(i)\stackrel{(\a^*)}\lto\bigoplus_{\a\colon j\to i}
P(j)\stackrel{(\a_*)}\lto C^-P(i)\lto 0 \quad\text{($i$ a sink)}\\
\label{eq:3}0\lto P(i)\xto{\smatrix{(\a^*)\\ (\b_*)}}\Big(\bigoplus_{\a\colon j\to i}
P(j)\Big)\oplus \Big(\bigoplus_{\b\colon i\to j}
C^-P(j)\Big)\xto{\smatrix{(\a_*)&(C^-\b^*)}} C^-P(i)\lto 0
\end{gather}
\end{lem}
\begin{proof}
For the first sequence, see Lemma~\ref{le:max}. 

The second sequence is
obtained from the first as follows. Let $i$ be a sink of $Q$. Then $i$
is a source of $\tilde Q=\s_iQ$. Now use Lemma~\ref{le:monoepi} and
apply $S_i^-$ to the sequence \eqref{eq:1} for $\tilde Q$ to obtain
the sequence \eqref{eq:2} for $Q$.

The third sequence is obtained from the second as follows.  Assume for
simplicity that $1,\ldots,n$ is an admissible labeling of the vertices
of $Q$. Then $i$ is a sink for $\tilde Q= \s_{i-1}\ldots\s_{1}Q$ and
we recall that $P(i)\cong S^-_1\ldots S^-_{i-1}S(i)$ by
Lemma~\ref{le:dimproj}.  Now use again Lemma~\ref{le:monoepi} and apply
$S^-_1\ldots S^-_{i-1}$ to the sequence \eqref{eq:2} for $\tilde Q$ to
obtain the sequence \eqref{eq:3} for $Q$.
\end{proof}

\begin{lem}\label{le:rel}
Suppose $Q$ has no oriented cycles. Given a vertex $i$ and $r\in\bbZ$, we have
\begin{equation}\label{eq:rel}
\sum_{\substack{\a\in Q_1\\ t(\a)=i}}C^r\a_*C^r\a^*+
\sum_{\substack{\a\in Q_1\\ s(\a)=i}}C^{r-1}\a^*C^r\a_*=0.
\end{equation}
\end{lem}
\begin{proof}
The exactness of the sequence \eqref{eq:3} implies
$\smatrix{(\a_*)&(C^-\b^*)}\smatrix{(\a^*)\\ (\b_*)}=0$ and therefore
\[\smatrix{(C^r\a_*)&(C^{r-1}\b^*)}\smatrix{(C^r\a^*)\\ (C^r\b_*)}=C^r\Big(\smatrix{(\a_*)&(C^-\b^*)}\smatrix{(\a^*)\\ (\b_*)}\Big)=0,\]
where $\a$ and $\b$ run through all arrows with $t(\a)=i=s(\b)$.  The
latter identity is precisely \eqref{eq:rel}.
\end{proof}

\subsection{Morphisms between preprojective representations}
We define a new quiver $\Ga=\bbZ Q$ as follows. Let
$$\Ga_0=\{i[r]\mid i\in Q_0,r\in\bbZ\}.$$ For each arrow $\a\colon
i\to j$ in $Q$ and $r\in\bbZ$, we have in $\Ga$ a pair of arrows
$$\a^*[r]\colon j[r]\lto i[r]\quad\text{and}\quad \a_*[r]\colon
i[r]\lto j[r-1].$$ Thus
$$\quad\Ga_1=\{\a^*[r]\mid\a\in
Q_1,r\in\bbZ\}\cup\{\a_*[r]\mid\a\in
Q_1,r\in\bbZ\}.$$

\setlength{\unitlength}{2.0pt}
\begin{figure}[!ht]
\begin{picture}(-20,50)

\put(-50,25){\circle*{1.2}}
\put(-30,25){\circle*{1.2}}
\put(-10,25){\circle*{1.2}}
\put(10,25){\circle*{1.2}}
\put(30,25){\circle*{1.2}}

\put(-60,05){\circle*{1.2}}
\put(-40,05){\circle*{1.2}}
\put(-20,05){\circle*{1.2}}
\put(00,05){\circle*{1.2}}
\put(20,05){\circle*{1.2}}
\put(40,05){\circle*{1.2}}

\put(-60,25){\circle*{1.2}}
\put(-40,25){\circle*{1.2}}
\put(-20,25){\circle*{1.2}}
\put(00,25){\circle*{1.2}}
\put(20,25){\circle*{1.2}}
\put(40,25){\circle*{1.2}}

\put(-60,45){\circle*{1.2}}
\put(-40,45){\circle*{1.2}}
\put(-20,45){\circle*{1.2}}
\put(00,45){\circle*{1.2}}
\put(20,45){\circle*{1.2}}
\put(40,45){\circle*{1.2}}

\put(-50,15){\circle*{1.2}}
\put(-30,15){\circle*{1.2}}
\put(-10,15){\circle*{1.2}}
\put(10,15){\circle*{1.2}}
\put(30,15){\circle*{1.2}}

\put(-50,35){\circle*{1.2}}
\put(-30,35){\circle*{1.2}}
\put(-10,35){\circle*{1.2}}
\put(10,35){\circle*{1.2}}
\put(30,35){\circle*{1.2}}

\put(01,06){\vector(1,1){8}}
\put(01,26){\vector(1,1){8}}
\put(01,06){\vector(1,1){8}}

\put(-59,24){\vector(1,-1){8}}
\put(-59,44){\vector(1,-1){8}}

\put(-39,24){\vector(1,-1){8}}
\put(-39,44){\vector(1,-1){8}}

\put(-19,24){\vector(1,-1){8}}
\put(-19,44){\vector(1,-1){8}}

\put(01,24){\vector(1,-1){8}}
\put(01,44){\vector(1,-1){8}}

\put(21,24){\vector(1,-1){8}}
\put(21,44){\vector(1,-1){8}}

\put(41,24){\vector(1,-1){8}}
\put(41,44){\vector(1,-1){8}}

\put(-69,14){\vector(1,-1){8}}
\put(-69,34){\vector(1,-1){8}}

\put(-49,14){\vector(1,-1){8}}
\put(-49,34){\vector(1,-1){8}}

\put(-29,14){\vector(1,-1){8}}
\put(-29,34){\vector(1,-1){8}}

\put(-09,14){\vector(1,-1){8}}
\put(-09,34){\vector(1,-1){8}}

\put(11,14){\vector(1,-1){8}}
\put(11,34){\vector(1,-1){8}}

\put(31,14){\vector(1,-1){8}}
\put(31,34){\vector(1,-1){8}}

\put(-69,16){\vector(1,1){8}}
\put(-69,36){\vector(1,1){8}}

\put(-49,16){\vector(1,1){8}}
\put(-49,36){\vector(1,1){8}}

\put(-29,16){\vector(1,1){8}}
\put(-29,36){\vector(1,1){8}}

\put(-09,16){\vector(1,1){8}}
\put(-09,36){\vector(1,1){8}}

\put(11,16){\vector(1,1){8}}
\put(11,36){\vector(1,1){8}}

\put(31,16){\vector(1,1){8}}
\put(31,36){\vector(1,1){8}}

\put(-59,06){\vector(1,1){8}}
\put(-59,26){\vector(1,1){8}}

\put(-39,06){\vector(1,1){8}}
\put(-39,26){\vector(1,1){8}}

\put(-19,06){\vector(1,1){8}}
\put(-19,26){\vector(1,1){8}}

\put(01,06){\vector(1,1){8}}
\put(01,26){\vector(1,1){8}}

\put(21,06){\vector(1,1){8}}
\put(21,26){\vector(1,1){8}}

\put(41,06){\vector(1,1){8}}
\put(41,26){\vector(1,1){8}}

\put(-69,25){\vector(1,0){8}}
\put(-59,25){\vector(1,0){8}}
\put(-49,25){\vector(1,0){8}}
\put(-39,25){\vector(1,0){8}}
\put(-29,25){\vector(1,0){8}}
\put(-19,25){\vector(1,0){8}}
\put(-09,25){\vector(1,0){8}}
\put(1,25){\vector(1,0){8}}
\put(11,25){\vector(1,0){8}}
\put(21,25){\vector(1,0){8}}
\put(31,25){\vector(1,0){8}}
\put(41,25){\vector(1,0){8}}

\put(-82,15){\circle*{.6}}
\put(-80,15){\circle*{.6}}
\put(-78,15){\circle*{.6}}

\put(-82,25){\circle*{.6}}
\put(-80,25){\circle*{.6}}
\put(-78,25){\circle*{.6}}

\put(-82,35){\circle*{.6}}
\put(-80,35){\circle*{.6}}
\put(-78,35){\circle*{.6}}

\put(62,15){\circle*{.6}}
\put(60,15){\circle*{.6}}
\put(58,15){\circle*{.6}}

\put(62,25){\circle*{.6}}
\put(60,25){\circle*{.6}}
\put(58,25){\circle*{.6}}

\put(62,35){\circle*{.6}}
\put(60,35){\circle*{.6}}
\put(58,35){\circle*{.6}}

\end{picture}
\caption*{\footnotesize Example: $\bbZ Q$ for $Q$ of Dynkin type $E_6$}
\end{figure}
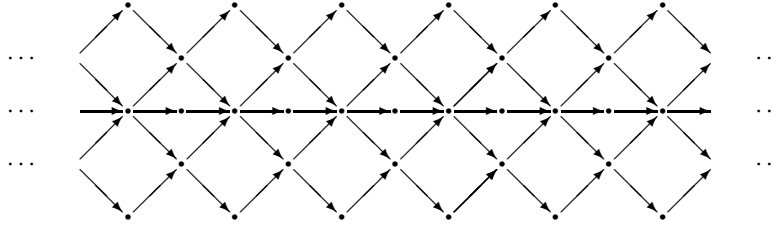

Observe that there exists a chain of irreducible morphisms
$$C^{r_1}P(i_1)\lto C^{r_2}P(i_2)\lto\cdots \lto C^{r_m}P(i_m)$$ if
and only if there is a path
$$i_1[r_1]\lto i_2[r_2]\lto \cdots\lto i_m[r_m]$$ in $\Ga$.  Here we
assume that $C^{r_j}P(i_j)\neq 0$ for $1\le j\le m$.  This observation
follows from Proposition~\ref{pr:irrpp} and motivates the following
construction.

We define for each pair of vertices $i[r]$ and $j[s]$ in $\Ga$ a
linear map
\begin{equation}\label{eq:map}
\pi\colon k[\Ga(i[r],j[s])]\lto\Hom(C^rP(i),C^sP(j))
\end{equation} 
by induction on the path length. For a path $\xi\colon i[r]\to j[s]$
in $\Ga$, let
$$\pi(\xi)=
\begin{cases}
\id_{C^rP(i)}&\text{if }\xi=\e_{i[r]},\\
C^r\a^*&\text{if }\xi=\a^*[r],\\
C^r\a_*&\text{if }\xi=\a_*[r],\\
\pi(\xi_l)\ldots\pi(\xi_1)&\text{if }\xi=\xi_l\ldots\xi_1, l>1.
\end{cases}$$

\begin{prop}\label{pr:pi}
Let $C^rP(i)$ and $C^sP(j)$ be two indecomposable preprojective
representations.  Then the linear map \eqref{eq:map} is an
epimorphism and its kernel is spanned by all elements of the
form
\begin{equation}\label{eq:ker}
\sum_{\substack{\a\in Q_1\\ t(\a)=l}}\t\a_*[t]\a^*[t]\s+
\sum_{\substack{\a\in Q_1\\ s(\a)=l}}\t\a^*[t-1]\a_*[t]\s,
\end{equation}
where $l[t]$ runs through all vertices of $\Ga$ and $\s,\t$ run through all paths
$\s\colon i[r]\to l[t]$ and $\t\colon l[t-1]\to j[s]$ in $\Ga$.
\end{prop}
\begin{proof}
Lemma~\ref{le:rel} implies that the kernel contains all elements of
the form \eqref{eq:ker}.  For the complete proof, see the discussion
of preprojective components in \cite[2.3.2]{R}.
\end{proof}

\subsection{The Dynkin case}
The following theorem summarises the structure of the morphisms
between representations of quivers of Dynkin type. Note that in this
case each indecomposable representation is preprojective and therefore
of the form $C^rP(i)$ for some vertex $i$ and some $r\le 0$.

\begin{thm}\label{th:mor}
Let $Q$ be a quiver whose underlying graph is a Dynkin
diagram. Suppose $X\cong C^rP(i)$ and $Y\cong C^sP(j)$ are two
indecomposable representations.
\begin{enumerate}
\item We have $\End(X)\cong k$ and every non-isomorphism $X\to Y$ is a
sum of compositions of irreducible morphisms between indecomposable
representations.
\item There exists an irreducible morphism $X\to Y$ if and only if
\begin{enumerate}
\item $r=s$ and there exists an arrow $j\to i$, or 
\item $r=s+1$ and there exists an arrow $i\to j$.
\end{enumerate}
Given two irreducible morphisms $\p,\p'\colon X\to Y$, there exists
$\la\in k$ with $\p'=\la\p$.
\item There exists an integer $d=d(X,Y)\ge 0$ such that 
\begin{equation*}\Hom(X,Y)=\Rad^0(X,Y)=\ldots=\Rad^d(X,Y)\supseteq \Rad^{d+1}(X,Y)=0
\end{equation*}
\end{enumerate}
\end{thm}
\begin{proof}
(1) We have $\End(X)\cong k$ by Proposition~\ref{pr:end}. There are
only finitely many isomorphism classes of indecomposable
representations by Theorem~\ref{th:Gabriel}. Thus the length of the
indecomposable representations is bounded by some $n$ and we have
therefore $\Rad^{2^n-1}(X,Y)=0$ by Proposition~\ref{pr:rad}. It
follows from Proposition~\ref{pr:irr} that every non-isomorphism $X\to
Y$ is a sum of compositions of irreducible morphisms between
indecomposable representations.

(2) See Proposition~\ref{pr:irrpp}. The last assertion follows from
    (3), because $\Irr(X,Y)\neq 0$ implies $\Rad^2(X,Y)= 0$.

(3) Let $d=d(i[r],j[s])$ be the length of a path $\xi\colon i[r]\to
j[s]$ in $\Ga$ and put $d=0$ if there is no such path. Observe that
$d$ does not depend on the choice of $\xi$ because any two parallel
paths in $\Ga$ have the same length. Here one uses that the underlying
graph of $Q$ is a tree. Now we apply (1). The assertion is clear if
$\Rad(X,Y)=0$.  A non-zero morphism $\p\in\Rad(X,Y)$ can be written as
a finite sum $\p=\sum_l\p_{ld}\ldots\p_{l2}\p_{l1}$ of compositions of
irreducible morphisms between indecomposable representations, and each
chain has length $d$ since it corresponds to a path $i[r]\to j[s]$
in $\Ga$.  Thus we have $\p\in\Rad^d(X,Y)$ but
$\p\not\in\Rad^{d+1}(X,Y)$.
\end{proof}

\section{The infinite radical}

In this section we characterise the quivers of finite representation
type in terms of morphisms between their representations. We use some
global properties, in particular infinite chains of radical morphisms.

Throughout this section we fix a quiver $Q$.

\subsection{Infinite chains of morphisms}

\begin{prop}
\label{pr:chain}
Let $Q$ be a quiver of Euclidean type. Then there exists an infinite
family of non-isomorphisms $\p_p\colon X_p\to X_{p+1}$, $p\ge 1$,
between indecomposable representations such that $\p_n\ldots\p_1\neq
0$ for all $n\ge 1$.
\end{prop}
\begin{proof}
Suppose first $Q$ is of type $\tilde A_n$. We allow any orientation,
in particular an oriented cycle.  We have an infinite family of
indecomposable representations $X(p)$, $p\ge 1$, by
Proposition~\ref{pr:cycle}. The construction of $X(p)$ shows that the
canonical inclusion $k^p\to k^{p+1}$ induces a monomorphism $\p_p\colon X(p)\to
X(p+1)$ for each $p\ge 1$.

Now suppose $Q$ has no oriented cycles.  We use for all $i\in Q_0$ the
monomorphism \[\m(i)\colon
P(i)\xto{\smatrix{(\a^*)\\ (\b_*)}}\Big(\bigoplus_{\a\colon j\to i}
P(j)\Big)\oplus \Big(\bigoplus_{\b\colon i\to j} C^-P(j)\Big)=E(i)\] in
\eqref{eq:3} and observe that $C^r\m(i)$ is a monomorphism for all
$r\le 0$. This follows from Lemma~\ref{le:monoepi} and the fact that
$C^rP(i)\neq 0$ for all $r\le 0$; see Theorem~\ref{th:tame}.  Now
choose a vertex $i_1$ of $Q$.  Let $X_1=P(i_1)$ and denote by
$\chi_1\colon X_1\to I$ a non-zero morphism to the indecomposable
injective representation $I=I(i_1)$. Then $\chi_1$ factors through
$\m(i_1)$ by Remark~\ref{re:inj}, because $\m(i_1)$ is a
monomorphism. Thus we can choose an indecomposable direct summand
$X_2=C^{r_2}P(i_2)$ of $E(i_1)$ corresponding to an arrow $\a_1$ and a
morphism $\chi_2\colon X_2\to I$ such that $\chi_2\p_1\neq 0$ where
$\p_1=\m(i_1)_{\a_1}$.  The morphism $\chi_2$ factors through
$C^{r_2}\m(i_2)$ because $C^{r_2}\m(i_2)$ is a monomorphism. Thus we
can choose an indecomposable direct summand $X_3=C^{r_3}P(i_3)$ of
$C^{r_2}E(i_2)$ corresponding to an arrow $\a_2$ and a morphism
$\chi_3\colon X_3\to I$ such that $\chi_3\p_2\p_1\neq 0$ where
$\p_2=C^{r_2}\m(i_2)_{\a_2}$. We continue and obtain an infinite
family of morphisms $\p_p\colon X_p\to X_{p+1}$ such that
$\p_n\ldots\p_1\neq 0$ for all $n\ge 1$.
\end{proof}

\subsection{A characterisation of finite representation type}
\begin{thm}
For a quiver $Q$, the following are equivalent:
\begin{enumerate}
\item The number of isomorphism classes of
indecomposable representations is finite.
\item There is a global bound for the length of every indecomposable
representation.
\item We have $\bigcap_{n\ge 0}\Rad^n(X,Y)=0$ for every pair $X,Y$ of representations.
\item Given an infinite family of non-isomorphism  $\p_i\colon X_i\to
X_{i+1}$, $i\ge 1$, between indecomposable representations, there exists
$n\ge 1$ such that $\p_n\ldots\p_1=0$.
\end{enumerate}
\end{thm}
\begin{proof}
(1) $\Rightarrow$ (2): Clear.

(2) $\Rightarrow$ (3): Use Proposition~\ref{pr:rad}. 

(3) $\Rightarrow$ (4): Suppose $\psi_n=\p_n\ldots\p_1\neq 0$ for all
$n\ge 1$. Then there exists $r\ge 1$ such that $\bigcup_{n\ge 1}\Ker
\psi_n=\Ker\psi_r \neq X_1$ because $X_1$ is finite dimensional. We
let $X=X_1/{\Ker\psi_r}$ and denote by $\p\colon X\to X_{r+1}$ the
canonical monomorphism. Then the composition $\p_n\ldots\p_{r+1}\p$ is
a monomorphism for all $n>r$. Choose an indecomposable injective
representations $I$ and a non-zero morphism $\chi\colon X\to I$ which
exist by Lemma~\ref{le:hom}. Then Remark~\ref{re:inj} implies that
$\chi$ factors through $\p_n\ldots\p_{r+1}\p$ for all $n>r$. Thus
$\chi$ belongs to $\bigcap_{n\ge 0}\Rad^n(X,I)$.

(4) $\Rightarrow$ (1): Suppose there are infinitely many
indecomposable representations. Then $Q$ contains an Euclidean
subquiver $Q'$ by Corollary~\ref{co:gabriel}. We obtain from
Proposition~\ref{pr:chain} an infinite chain $\p_i\colon X_i\to
X_{i+1}$ of non-isomorphisms between indecomposable representations of
$Q'$, which we can extend to a chain of morphisms for $Q$ as in the
proof of Corollary~\ref{co:gabriel}.
\end{proof}

\section{Regular representations}

We study regular representations and concentrate on two particular
quivers: the Jordan quiver and the Kronecker quiver.  In both cases we
provide complete classifications of all representations.  Note that
the Kronecker quiver is the easiest but also most important case among
the quivers of Euclidean type without oriented cycles.

\subsection{Direction of morphisms}

Let $Q$ be a quiver without oriented cycles. We call a representation
\emph{regular} if all its indecomposable direct summands are regular.

\begin{lem}\label{le:hom-reg1}
Let $X,Y$ be indecomposable representations.
\begin{enumerate}
\item If $Y$ is (pre)projective and $X$ is not, then $\Hom(X,Y)=0$.
\item If $Y$ is (pre)injective and $X$ is not, then $\Hom(Y,X)=0$.
\end{enumerate}
\end{lem}

\begin{proof}
(1) A preprojective representation $Y$ is of the form $C^r P(i)$ for
  some $r\le 0$ and some $i\in Q_0$.  We fix a morphism $\p\colon X\to
  Y$ and need to show that $\p=0$.

Consider first the case $r=0$, that is, $Y$ is projective and $X$
is not. We use induction on the length $l=\ell(i)$ of the longest
path in $Q$ starting at $i$. Note that $\p$ is not an epimorphism,
since every epimorphism splits by Remark~\ref{re:inj}. If $l=0$, then
$P(i)$ is simple and therefore $\p=0$. If $l>0$, then $\p$ factors
through $\bigoplus_{\a\colon i\to j}P(j)\xto{(\a^*)} P(i)$ by
Lemma~\ref{le:max}. We have $\ell(j)<l$ for each arrow $i\to j$, and
therefore $\p=0$.

Now suppose that $r<0$ and that $X$ is not preprojective. Thus we have
$X\cong C^rC^{-r} X$ by Proposition~\ref{pr:proj}, and therefore
\[\Hom(X,Y)\cong\Hom(C^rC^{-r}X,Y)\cong\Hom(C^{-r}X,P(i))=0,\]
where the second isomorphism follows from Proposition~\ref{pr:coxmor}
and the last identity follows from the first part of this proof. Thus
$\p=0$.

(2) is dual to (1). 
\end{proof}

\begin{rem}
Let $X,Y$ be indecomposable representations such that $X=C^rP(i)$ is preprojective and $Y$ is not. Then

\[\Hom(C^rP(i),Y)\cong\Hom(P(i),C^{-r}Y)\cong (C^{-r}Y)_i.\]
\end{rem}

\begin{lem}\label{le:hom-reg2}
Let $Q$ be a quiver of Euclidean type. For any morphism $\p\colon X\to
Y$ between regular representations, the representations $\Ker\p$,
$\Coker\p$, and $\Im\p$ are regular.
\end{lem}
\begin{proof}
We apply Lemma~\ref{le:hom-reg1}. 
The image $\Im \p$ is a subrepresentation of $Y$ and admits therefore
no preinjective direct summand. On the other hand, $\Im\p$ is a
quotient of $X$ and has therefore no preprojective direct
summand. Thus $\Im\p$ is regular.

In order to show that $\Ker\p$ is regular, we compute the defect and
apply Proposition~\ref{pr:defect}.  We have $\partial\Ker\p=\partial
X-\partial \Im\p=0$ since $\Im\p$ is regular.  $\Ker \p$ is a
subrepresentation of $X$ and admits therefore no preinjective direct
summand. Thus any preprojective direct summand $P$ of $\Ker\p$ would
imply $\partial\Ker\p\le\partial P<0$. It follows that $\Ker\p$ is
regular, and the dual argument shows that $\Coker\p$ is regular.
\end{proof}

\subsection{Jordan quiver representations}

We study the representations of the following Jordan quiver
\[\xymatrix{\circ\ar@(ur,dr)}\]
and assume throughout that the field $k$ is algebraically closed. By
definition, such a representation is a pair $(V,\p)$ consisting of a
vector space $V$ and an endomorphism $\p\colon V\to V$. For each
integer $p\ge 1$ and $\la\in k$, we have the representation
\[J_{p,\la}=(k^p,J(p,\la))\quad \text{where} \quad
J(p,\la)=
\smatrix{\la&1&0&\cdots&0\\0&\la&1&\cdots&0\\ 
\vdots&\vdots&\vdots&&\vdots\\0&0&0&\cdots&1\\0&0&0&\cdots&\la}\]
denotes the Jordan block of size $p$ with eigenvalue $\la$ (sending a
standard basis vector $e_i$ to $\la e_i +e_{i-1}$).

\begin{thm}[Jordan normal form]\label{th:jordan}
\pushQED{\qed}
The finite dimensional indecomposable representations of the Jordan
quiver are, up to isomorphism, precisely the representations  $J_{p,\la}$
with $p\ge 1$ and $\la\in k$.\qedhere
\end{thm}

Fix $\la\in k$. For each pair of integers $p,q\ge 1$, we define a
\emph{standard morphism}
\[\p_{p,q}\colon J_{p,\la}\lto J_{q,\la}\quad\text{with}\quad\p_{p,q}(e_i)=
\begin{cases}e_i&\text{if }p\le q,\\e_{i-(p-q)}&\text{if }p>q,\end{cases}\]
where $\{e_1,\ldots,e_p\}$ denotes the standard basis of $k^p$ and
$e_i=0$ for $i\le 0$. If $p>q$, these morphisms induce an exact
sequence
\[0\lto J_{p-q,\la}\xto{\p_{p-q,p}} J_{p,\la}\xto{\p_{p,q}}J_{q,\la}\lto 0.\]

Let us identify for each $q<p$ the representation $J_{q,\la}$ with the
image of $\p_{q,p}$. Thus we obtain a chain of subrepresentations
\[J_{1,\la}\subseteq J_{2,\la}\subseteq\ldots\subseteq J_{p-1,\la} 
\subseteq J_{p,\la}.\]

A representation is called \emph{uniserial} if its subrepresentations
are linearly ordered, that is, for each pair of
subrepresentations $U,V$ one has $U\subseteq V$ or $V\subseteq U$.

\begin{lem}\label{le:uniserial}
Each representation $J_{p,\la}$ is uniserial. More precisely, each
proper subrepresentation is of the form $J_{q,\la}$ for some $1\le q< p$, and
$J_{p,\la}/J_{q,\la}\cong J_{p-q,\la}$.
\end{lem}
\begin{proof}
The proof is by induction on $p$. The assertion is clear for $p=1$; so
let $p>1$.  Each proper $J(p,\la)$-invariant subspace of $k^p$
contains $e_1$. Thus each proper subrepresentation of $J_{p,\la}$
contains $J_{1,\la}$, which is the kernel of $\p_{p,p-1}$. The
representation $J_{p-1,\la}$ is uniserial and the epimorphism
$\p_{p,p-1}$ induces an inclusion preserving bijection between the
subrepresentations of $J_{p,\la}$ containing $J_{1,\la}$ and all
subrepresentations of $J_{p-1,\la}$. This yields the description of
all subrepresentations of $J_{p,\la}$. The isomorphism
$J_{p,\la}/J_{q,\la}\xto{\sim} J_{p-q,\la}$ is induced by
$\p_{p,p-q}$.
\end{proof}

\begin{lem}\label{le:hom-jordan}
Fix two representations $J_{p,\la}$ and $J_{q,\mu}$.
\begin{enumerate}
\item Let $\la=\mu$. Then $\{\p_{i,q}\p_{p,i}\mid 1\le i\le \min(p,q)\}$ is a basis
  for $\Hom(J_{p,\la},J_{q,\mu})$.
\item Let $\la\neq\mu$. Then $\Hom(J_{p,\la},J_{q,\mu})=0$ and $\Ext(J_{p,\la},J_{q,\mu})=0$.
\end{enumerate}
\end{lem}
\begin{proof}
(1) Let $\la=\mu$ and set $\Phi_{p,q}=\{\p_{i,q}\p_{p,i}\mid 1\le i\le
  \min(p,q)\}$. The proof is by induction on $q$. The assertion is
  clear for $q=1$; so let $q>1$.  The exact sequence
\[0\lto J_{q-1,\mu}\xto{\p_{q-1,q}} J_{q,\mu}\xto{\p_{q,1}}J_{1,\mu}\lto 0\]
induces an exact sequence
\begin{equation}\label{eq:hom}
0\lto \Hom(J_{p,\la},J_{q-1,\mu})\xto{(J_{p,\la},\p_{q-1,q})} \Hom(J_{p,\la},J_{q,\mu})
\xto{(J_{p,\la},\p_{q,1})}\Hom(J_{p,\la},J_{1,\mu})
\end{equation} 
by Lemma~\ref{le:hom-seq}. The map $(J_{p,\la},\p_{q-1,q})$ takes the
basis $\Phi_{p,q-1}$ of $\Hom(J_{p,\la},J_{q-1,\mu})$ to
$\Phi=\{\p_{i,q}\p_{p,i}\mid 1\le i\le \min(p,q-1)\}$ which is a
subset of $\Phi_{p,q}$. Using the the exactness of \eqref{eq:hom},
it follows that $\Phi_{p,q}=\Phi\cup\{\p_{p,q}\}$ is a
basis of $\Hom(J_{p,\la},J_{q,\mu})$.

(2) The proof goes by induction on $p$ and $q$, using for $\Hom(-,-)$
the sequence \eqref{eq:hom} and its contravariant counterpart. For
$\Ext(-,-)$ one uses Lemma~\ref{le:ext}. The case $p=1=q$ is clear.
%
\end{proof}

\begin{lem}
For a morphism $\p\colon J_{p,\la}\to J_{q,\la}$, the following are
equivalent:
\begin{enumerate}
\item The morphism $\p$ is irreducible.
\item The representation $\Ker\p\oplus\Coker\p$ is simple.
\item $|p-q|=1$ and $\p$ is a monomorphism or an epimorphism.
\end{enumerate}
\end{lem}
\begin{proof}
(1) $\Rightarrow$ (2): An irreducible morphism is either a
  monomorphism or an epimorphism. It suffices to discuss the case that
  $\p$ is an epimorphism; the other case is dual. If $\Ker\p$ is not
  simple and $S\subseteq\Ker \p$ is a simple subrepresentation, then $\p$ can
  be written as composite $J_{p,\la}\to J_{p,\la}/S\to J_{q,\la}$ of two proper
  epimorphisms. This is a contradiction, and therefore
  $\Ker\p\oplus\Coker\p$ is simple.

(2) $\Rightarrow$ (3): Clear.

(3) $\Rightarrow$ (1): It suffices to consider an epimorphism
  $J_{q+1,\la}\to J_{q,\la}$; the dual argument works for a
  monomorphism $J_{p,\la}\to J_{p+1,\la}$. Let $J_{q+1,\la}\xto{\a}
  X\xto{\b}J_{q,\la}$ be a factorisation and fix a decomposition
  $X=\bigoplus_iX_i$ into indecomposable representations. Then
  $\b_{i_0}\a_{i_0}$ is an epimorphism for at least one index
  $i_0$. It follows that
  $X_{i_0}=J_{r,\la}$ for some $r\ge q$. If $r=q$, then $\b_{i_0}$ is
  an isomorphism, and therefore $\b$ is a split
  epimorphism. Otherwise, we obtain a factorisation
  $J_{q+1,\la}\xto{\a_{i_0}}
  X_{i_0}\xto{\b'_{i_0}}J_{q+1,\la}\xto{\p_{q+1,q}} J_{q,\la}$ of the
  epimorphism $\b_{i_0}\a_{i_0}$. It follows that $\b'_{i_0}\a_{i_0}$
  is an epimorphism and hence an isomorphism. Thus $\a$ is a split
  monomorphism.
\end{proof}

\subsection{Kronecker quiver representations}

We study the representations of the following Kronecker quiver
\[\xymatrix{ 1\ar@<3.0pt>[rr]\ar@<-3.0pt>[rr]&& 2}\]
and assume throughout that the field $k$ is algebraically closed. 
By definition, such a representation is a pair of linear maps between
two vector spaces. Let us list the indecomposable representations.

For each integer $r\ge 0$, there are the preprojective and preinjective
representations:
\[\xymatrix{P_r:\quad k^r\ar@<3.0pt>[rr]^-{\smatrix{\id\\0}}
\ar@<-3.0pt>[rr]_-{\smatrix{0\\ \id}}&&k^{r+1}}\qquad
\xymatrix{I_r:\quad k^{r+1}\ar@<3.0pt>[rr]^-{\smatrix{\id&0}}
  \ar@<-3.0pt>[rr]_-{\smatrix{0&\id}}&&k^{r}}\] where in each block
matrix $\id=\id_{k^r}$. We have $C^{-r}P(1)=P_{2r+1}$ and
$C^{-r}P(2)=P_{2r}$, while $C^{r}I(1)=I_{2r}$ and
$C^{r}I(2)=I_{2r+1}$; see Theorem~\ref{th:tame}.

The regular representations are indexed by points of the
\emph{projective line} $\bbP^1(k)$, which are by definition pairs
$(\la_0:\la_1)$ of elements in $k$ different from $(0:0)$ and subject
to the relation $(\la_0:\la_1)=(\alpha \la_0:\alpha \la_1)$ for all
$\alpha\in k$, $\alpha\neq 0$.

For each integer $p\ge 1$ and $\la\in\bbP^1(k)$,
we consider the representations
\[\xymatrix{k^p\ar@<3.0pt>[rr]^-{J(p,\la_0)}\ar@<-3.0pt>[rr]_-{\id}&&k^p\qquad
\text{for $\la=(\la_0:1)$, and}\qquad
k^p\ar@<3.0pt>[rr]^-{\id}\ar@<-3.0pt>[rr]_-{J(p,\la_1)}&&k^p\qquad\text{for
  $\la=(1:\la_1)$,}}\] where $J(p,\mu)$ denotes the Jordan block of
size $p$ with eigenvalue $\mu$.  Observe that both representations are
isomorphic in case $\la_0\neq 0\neq\la_1$.  To see this, note that the
first representation is isomorphic
to \[\xymatrix{k^p\ar@<3.0pt>[rr]^-{\id}\ar@<-3.0pt>[rr]_-{J(p,\la_0)^{-1}}&&k^p},\]
and that $J(p,\la_0)^{-1}\cong J(p,\la_1)$ since both endomorphisms
are indecomposable with same eigenvalue.  Thus we denote this
representation by $R_{p,\la}$, and it is well-defined up to an
isomorphism.

\begin{thm}[Kronecker]\label{th:Kron}
The representations $P_r$, $I_r$ ($r\ge 0$), and $R_{p,\la}$ ($p\ge 1$
and $\la\in\bbP^1(k)$) form, up to isomorphism, a complete list of
finite dimensional indecomposable representations of the Kronecker
quiver.
\end{thm}
\begin{proof}
The preprojective and preinjective representations have been
classified in Theorem~\ref{th:tame}, while the regular ones are
classified further below in Proposition~\ref{pr:reg-Kro}.
\end{proof}

Fix a point $\la\in\bbP^1(k)$ with $\la=(\la_0:1)$ or
$\la=(1:\la_1)$. Then for each pair $p,q\ge 1$ of integers, we have an
obvious bijection
\[\Hom(J_{p,\la_i},J_{q,\la_i})\stackrel{\sim}\lto\Hom(R_{p,\la},R_{q,\la}).\]
In particular, a standard morphism $\p_{p,q}\colon J_{p,\la_i}\to
J_{q,\la_i}$ induces another standard morphism $R_{p,\la}\to
R_{q,\la}$ which we denote by $\psi_{p,q}$.

\begin{lem}\label{le:hom-kron}
Fix two representations $R_{p,\la}$ and $R_{q,\mu}$.
\begin{enumerate}
\item Let $\la=\mu$. Then $\{\psi_{i,q}\psi_{p,i}\mid 1\le i\le \min(p,q)\}$ is a basis
  for $\Hom(R_{p,\la},R_{q,\mu})$.
\item Let $\la\neq\mu$. Then $\Hom(R_{p,\la},R_{q,\mu})=0$ and $\Ext(R_{p,\la},R_{q,\mu})=0$.
\end{enumerate}
\end{lem}
\begin{proof}
The case $\la=\mu$ is clear from Lemma~\ref{le:hom-jordan}. For
$\la\neq\mu$, the same induction argument as in
Lemma~\ref{le:hom-jordan} can be used. It remains to check the case
$p=1=q$. It is easily checked that $\Hom(R_{1,\la},R_{1,\mu})=0$.
Thus we consider an exact sequence $0\to R_{1,\mu}\to E\to
R_{1,\la}\to 0$ and need to show that it splits. Suppose first that
$\la_1\neq 0\neq\mu_1$. We may assume that $\la_1=1=\mu_1$ and obtain
the following commutative diagram of linear maps with exact rows.
\[\xymatrix{ 0\ar[r]&k\ar[r]\ar@<-3pt>[d]_{\mu_0}\ar@<3pt>[d]^{\id}&
k^2\ar[r]\ar@<-3pt>[d]_{\p_0}\ar@<3pt>[d]^{\p_1}&
k\ar[r]\ar@<-3pt>[d]_{\la_0}\ar@<3pt>[d]^{\id}&0\\ 0\ar[r]&k\ar[r]&k^2\ar[r]&k\ar[r]&0}\]
It follows from Lemma~\ref{le:snake} that $\p_1$ is an isomorphism,
and we may assume that $\p_1=\id$. This yields an exact sequence $0\to
J_{1,\mu_0}\to E_0\to J_{1,\la_0}\to 0$ of Jordan quiver
representations which splits by Lemma~\ref{le:hom-jordan}, since
$\la_0\neq\mu_0$. Thus $\eta$ splits. The case $\la_0\neq 0\neq\mu_0$
is analogous, and the remaining case $\{\la,\mu\}=\{(1:0),(0:1)\}$ is
easy.
\end{proof}

\begin{lem}\label{le:reg-sub}
Any indecomposable regular representation admits a subrepresentation
that is isomorphic to  $R_{1,\la}$ for some $\la\in\bbP^1(k)$.
\end{lem}
\begin{proof}
Suppose the representation $X$ is given by a pair of linear maps
$\p,\psi\colon V\to W$.
Observe that $\dim V=\dim W$ since $X$ is regular. If $\p$ is
bijective, then we may assume it is the identity so that $\psi$ is a
Jordan block, say $\psi=J(p,\la_1)$. This implies $X\cong R_{p,\la}$,
where $\la=(1:\la_1)$. Thus $X$ admits  $R_{1,\la}$ as a subrepresentation.

Now suppose that $\p$ is not bijective and put $K=\Ker\p\cap\Ker\psi$.
Then
$K\,\genfrac{}{}{0pt}{}{\raisebox{-1.75pt}{$\lto$}}{\raisebox{1.75pt}{$\lto$}}\,0$
is a subrepresentation of $X$ and the inclusion is a split
monomorphism. Thus $K=0$, and therefore $X$ admits $R_{1,\la}$ with
$\la=(0:1)$ as a subrepresentation.
\end{proof}

\begin{prop}\label{pr:reg-Kro}
The finite dimensional indecomposable and regular representations of
the Kronecker quiver are, up to isomorphism, precisely the
representations $R_{p,\la}$ with $p\ge 1$ and $\la\in \bbP^1(k)$.
\end{prop}
\begin{proof}
Fix an indecomposable regular representation $X$ that is given by a
pair of linear maps $\p,\psi\colon V\to W$. The proof goes by
induction on the dimension of $X$. Lemma~\ref{le:reg-sub} yields an
exact sequence
\begin{equation}\label{eq:reg}
0\lto R_{1,\la}\lto X\lto Y\lto 0
\end{equation} 
for some $\la\in\bbP^1(k)$.  If $Y=0$, then $X\cong R_{1,\la}$. Thus
we assume $Y\neq 0$ and fix a decomposition $Y=\bigoplus_{i=1}^r Y_i$
into indecomposable representations. Note that $Y$ is regular by
Lemma~\ref{le:hom-reg2}, so that each $Y_i$ is of the form
$R_{q_i,\mu_i}$ for some pair $q_i,\mu_i$. The sequence \eqref{eq:reg}
induces for each $i$ an exact sequence
\[0\lto R_{1,\la}\lto X_i\lto R_{q_i,\mu_i}\lto 0,\]
where $X_i$ denotes the preimage of $Y_i$ under the epimorphism $X\to
Y$.  Note that this sequence does not split since $X$ is
indecomposable.  Thus $\mu_i=\la$ for all $i$ by
Lemma~\ref{le:hom-kron}. We may assume $\la=(\la_0:1)$ or
$\la=(1:\la_1)$, and therefore one of $\p$ and $\psi$ can be chosen to
be the identity, using Lemma~\ref{le:snake} as in the proof of
Lemma~\ref{le:hom-kron}. It follows that the other map is a Jordan
block, say $J(p,\la_i)$, and therefore $X\cong R_{p,\la}$.
\end{proof}

\section{Wild phenomena}

We consider the following $n$-subspace quiver:
\[\xymatrix@!0 @R=1.2em @C=2.4em {
&&&0\\ \La_n\qquad\\&1\ar[uurr]&2\ar[uur]&\cdots&n\ar[uul]}\] The
representations of $\La_n$ are basically configurations of $n$
subspaces of a fixed vector space. For $n\le 3$, the underlying diagram
is of Dynkin type, while it is of Euclidean type for $n=4$. In this
section we demonstrate some wild phenomenon for $n\ge 5$. The same
phenomenon occurs for the $r$-Kronecker quiver
\[\xymatrix{K_r&1\ar @/u2.75ex/@<3.0pt>[rr]^{\a_1}
\ar @/d2.75ex/@<-3.0pt>[rr]^{\vdots}_{\a_r}\ar[rr]^{\a_2}&& 2}\]
in case $r\ge 3$.

Throughout this section we fix a quiver $Q$.

\subsection{Total representations}

We assign to each representation $X$ of $Q$ its \emph{total
  representation}, that is, a single vector space together with a
distinguished set of endomorphisms. Put $\bar X=\bigoplus_{i\in
  Q_0}X_i$ and denote for each $i\in Q_0$ by $\bar X_i$ the canonical
projection $\bar X\twoheadrightarrow X_i\rightarrowtail \bar X$. For
$\a\in Q_1$ let $\bar X_\a$ denote the linear map $\bar
X\twoheadrightarrow X_{s(\a)}\xto{X_\a} X_{t(\a)}\rightarrowtail \bar
X$.  For a morphism $\p\colon X\to Y$ of representations, let
$\bar\p=(\p_i)_{i\in Q_0}$ be the induced linear map $\bar X\to\bar
Y$.

\begin{lem}
Sending a representation to its total reprsentation induces a
bijection between all representations of $Q$ and all families
$(V,\p_i,\p_\a)_{i\in Q_0,\a\in Q_1}$ of vector spaces $V$ with
endomorphisms $\p_i,\p_\a\colon V\to V$ satisfying
\[\sum_{i\in Q_0}\p_i=\id_V, \quad\p_i\p_j=\d_{ij}\p_i\quad (i,j\in Q_0),\quad\text{and}\quad
\p_{t(\a)}\p_{\a}\p_{s(\a)}=\p_{\a}\quad (\a\in Q_1). \]
\end{lem}
\begin{proof}
The inverse map sends a family $(V,\p_i,\p_\a)_{i\in Q_0,\a\in Q_1}$
to the representation $X$ with $X_i=\Im\p_i$ for $i\in Q_0$ and
$X_\a=(X_{s(\a)} \rightarrowtail V\xto{\p_\a} V\twoheadrightarrow
X_{t(\a)})$ for $\a\in Q_1$.
\end{proof}

\begin{lem}\label{le:hom-bar}
The map $\Hom(X,Y)\to \Hom(\bar X,\bar Y)$ sending a morphism $\p$ to
$\bar\p$ identifies $\Hom(X,Y)$ with the subspace of
all linear maps $\psi\colon\bar X\to\bar Y$ satisfying 
\[\psi\bar X_{i}=\bar Y_{i}\psi\quad (i\in Q_0)\quad\text{and}\quad \psi
\bar X_{\a}=\bar Y_{\a}\psi\quad(\a\in Q_1).\]
\end{lem}
\begin{proof} The inverse map sends $\psi\colon\bar X\to\bar Y$ to
the family $(X_i\rightarrowtail\bar X\xto{\psi}\bar
Y\twoheadrightarrow Y_i)_{i\in Q_0}$.
\end{proof}
\subsection{A representation embedding}

A functor $F\colon\Rep(Q,k)\to\Rep(Q',k)$ is called \emph{embedding}
if $F$ induces  a bijection
$\Hom(X,Y)\xto{\sim}\Hom(FX,FY)$ for each pair of representations $X,Y$ of $Q$.

We construct an embedding $F_Q\colon\Rep(Q,k)\to\Rep(K_3,k)$ as
composite of two embeddings
\[\Rep(Q,k)\stackrel{E}\lto\Rep(\Ga,k)\stackrel{F}\lto\Rep(K_3,k)\]
with intermediate quiver 
\[\xymatrix{\Ga\qquad&\circ\ar@(ur,dr)^\t\ar@(ul,dl)_\s}.\]
The embedding $F_Q$ shows that the representations of $K_3$ are as
complicated as the representations of any other quiver; in this sense
$\Rep(K_3,k)$ is \emph{wild}.

Let $Q_0=\{1,\ldots,n\}$ and $Q_1=\{\a_1,\ldots,\a_r\}$. Fix
representations $X,Y$ of $Q$ and define the representation $E X$ of $\Ga$
as follows. The underlying vector space is $\bar X^{n+r+2}$. Any endomorphism of
this space can be written as $(n+r+2)\times(n+r+2)$ block matrix with
each entry an element of $\End\bar X$. So we define
 \[(EX)_\s=\smatrix{
0&\id&0&\cdots&0\\
0&0&\id&\cdots&0\\
\vdots&\vdots&\vdots&&\vdots\\
0&0&0&\cdots&\id\\
0&0&0&\cdots&0
}\quad\text{and}\quad
(EX)_\t=\smatrix{
0&0&0&\cdots&0&0&0&\cdots&0&0&0\\
\id&0&0&\cdots&0&0&0&\cdots&0&0&0\\
\bar X_{1}&\id&0&\cdots&0&0&0&\cdots&0&0&0\\
0&\bar X_{2}&\id&\cdots&0&0&0&\cdots&0&0&0\\
\vdots&\vdots&\vdots&&\vdots&\vdots&\vdots&&\vdots&\vdots&\vdots\\
0&0&0&\cdots&\bar X_{n}&\id&0&\cdots&0&0&0\\
0&0&0&\cdots&0&\bar X_{\a_1}&\id&\cdots&0&0&0\\
0&0&0&\cdots&0&0&\bar X_{\a_2}&\cdots&0&0&0\\
\vdots&\vdots&\vdots&&\vdots&\vdots&\vdots&&\vdots&\vdots&\vdots\\
0&0&0&\cdots&0&0&0&\cdots&\bar X_{\a_r}&\id&0
}.
\]
For a morphism $\p\colon X\to Y$, define a morphism $EX\to EY$ by
\[E\p=\smatrix{\bar\p&0&\cdots&0\\0&\bar\p&\cdots&0\\ \vdots&\vdots&&\vdots\\ 
0&0&\cdots&\bar\p}.\]
 
\begin{lem}\label{le:hom-E}
$E$ induces an isomorphism $\Hom(X,Y)\xto{\sim}\Hom(EX,EY)$.
\end{lem}
\begin{proof}
The map is certainly injective. In order to show that it is
surjective, fix a linear map $\psi\colon \bar X^{n+r+2}\to \bar
Y^{n+r+2}$ and suppose it is a morphism of $\Ga$-representations. The
condition $\psi(EX)_\s=(EY)_\s\psi$ implies that $\psi$ is a
$(n+r+2)\times(n+r+2)$ block matrix of the
form \[\smatrix{\chi&0&\cdots&0\\0&\chi&\cdots&0\\ \vdots&\vdots&&\vdots\\ 0&0&\cdots&\chi}\]
for some linear map $\chi\colon\bar X\to\bar Y$.  The second condition
$\psi(EX)_\t=(EY)_\t\psi$ implies that $\chi \bar X_{i}=\bar
Y_{i}\chi$ for all $i\in Q_0$ and $\chi \bar X_{\a}=\bar Y_{\a}\chi$
for all $\a\in Q_1$. Thus $\chi=\bar\p$ for some morphism $\p\colon
X\to Y$ by Lemma~\ref{le:hom-bar}, and therefore $\psi=E\p$.
\end{proof}

Let $X,Y$ be representations of $\Ga$. We identify $X$ with its
underlying vector space and define a representation $F X$ of the quiver $K_3$ by
\[\xymatrix{X\ar @/u2.75ex/@<3.0pt>[rr]^{X_\s}
\ar @/d2.75ex/@<-3.0pt>[rr]^{\id}\ar[rr]^{X_\t}&& X}.\]
Given a morphism $\p\colon X\to Y$, define $F\p\colon FX\to FY$ by $(F\p)_1=\p=(F\p)_2$.

\begin{lem}\label{le:hom-F}
$F$ induces an isomorphism $\Hom(X,Y)\xto{\sim}\Hom(FX,FY)$.
\end{lem}
\begin{proof}
The inverse map sends a morphism $\p\colon FX\to FY$ to $\p_1\colon X\to Y$.
\end{proof}

\begin{prop}\label{pr:wild_K3}
$F_Q=FE$ yields an embedding  $\Rep(Q,k)\to\Rep(K_3,k)$.
\end{prop}
\begin{proof}
Combine Lemmas~\ref{le:hom-E} and \ref{le:hom-F}.
\end{proof}

It is clear that there is a similar embedding $\Rep(Q,k)\to\Rep(K_r,k)$ for each $r>3$.

\subsection{From Kronecker to subspace representations}
We construct for $r\ge 1$ an embedding
\[F_r\colon \Rep(K_r,k)\lto\Rep(\La_{r+2},k).\]
Let $X,Y$ be representations of the $r$-Kronecker quiver $K_r$.
Define $F_rX$ by specifying an $r+2$-subspace system as follows:
\begin{gather*}(F_rX)_0=X_1\times X_2,\quad (F_rX)_1=X_1\times\{0\},\quad (F_rX)_2=\{0\}\times X_2,\\
(F_rX)_{i+2}=\{(v,X_{\a_{i}}(v))\in X_1\times X_2\mid v\in X_1\} \quad (1\le i\le r).
\end{gather*}
Given a morphism $\p\colon X\to Y$, define $F_r\p\colon F_rX\to F_rY$ by
$(F_r\p)_i(v,w)=(\p_1 (v),\p_2 (w))$ for $0\le i\le r+2$.

\begin{lem}
$F_r$ induces an isomorphism $\Hom(X,Y)\xto{\sim}\Hom(F_rX,F_rY)$.
\end{lem}
\begin{proof}
The inverse map sends a morphism $\p\colon F_rX\to F_rY$ to 
a morphism $\psi\colon X\to Y$ given by $\psi_1=\p_1$ and $\psi_2=\p_2$.
\end{proof}

\begin{prop}
\pushQED{\qed}
$F_3F_Q$ yields an embedding  $\Rep(Q,k)\to\Rep(\La_5,k)$.\qedhere
\end{prop}

\section{Radical square zero representations}

We fix a quiver $Q$ and construct functors between radical square zero
representations of $Q$ and representations of the corresponding
separated quiver $Q^s$.

\subsection{The radical}

Let $X$ be a representation of $Q$. The \emph{radical} of $X$ is a
subrepresentation which we denote by $\Rad X$. We define $(\Rad
X)_i=\sum_{\a\colon j\to i}\Im X_\a$ for each vertex $i$.  The
\emph{Jacobson radical} $\rad X$ of $X$ is by definition the
intersection of all maximal subrepresentations of $X$.  For $n\ge 1$,
let $\Rad X^{n+1}=\Rad (\Rad^n X)$ and $\rad X^{n+1}=\rad (\rad^n X)$.
Note that each morphism $X\to Y$ induces morphisms $\Rad X\to \Rad Y$
and $\rad X\to \rad Y$.

\begin{lem}\label{le:jrad}
$\rad X\subseteq \Rad X$ and equality holds if $\Rad^n X=0$ for some $n\ge 1$.
\end{lem}
\begin{proof}
The representation $X/\Rad X$ is a direct sum of simple
representations, and therefore $\rad(X/\Rad X)=0$. The canonical
morphism $\pi\colon X\to X/\Rad X$ sends $\rad X$ to $\rad(X/\Rad
X)$. Thus $\rad X\subseteq \Rad X$.

Now suppose that $\Rad^n X=0$ for some $n\ge 1$.  Given a maximal
subrepresentation $U\subseteq X$, it follows that $\Rad (X/U)=0$. Thus
$\Rad X\subseteq U$, and therefore $\Rad X\subseteq\rad X$.
\end{proof}

Let $Q$ be a quiver without oriented cycles.  If $n$ equals the length
of the longest path in $Q$, then $\Rad^{n+1} X=0$ for each
representation $X$. On the other hand, any quiver with an oriented
cycle has a simple representation $X$ such that $\Rad
X=X$, while $\rad X=0$.

\subsection{The separated quiver}
We define the \emph{separated quiver} $Q^s$ as follows.  Let
$Q_0=\{1,2,\ldots,n\}$ be the vertices of $Q$ and set
$Q^s_0=\{1,\ldots,n,1',\ldots,n'\}$. For each arrow $\a\colon i\to j$
of $Q$, there is an arrow $\bar\a\colon i\to j'$ in $Q^s$. Thus
$Q_1^s=\{\bar\a\mid\a\in Q_1\}$.

There are two functors
\[S\colon\Rep(Q,k)\lto\Rep(Q^s,k)\quad
\text{and}\quad T\colon\Rep(Q^s,k)\lto\Rep(Q,k)\] which are defined as
follows. Given a representation $X$ of $Q$, let $(SX)_i=(X/\Rad X)_i$
and $(SX)_{i'}=(\Rad X/\Rad^2 X)_i$ for each vertex $i\in Q_0$. For
each arrow $\a\colon i\to j$ of $Q$, let
$(SX)_{\bar\a}\colon(SX)_i\to(SX)_{j'}$ be the map which is induced by
$X_\a$.  Given a representation $Y$ of $Q^s$, let $(TY)_i=Y_i\oplus
Y_{i'}$ for each vertex $i\in Q_0$. For each arrow $\a\in Q_1$, let
$(TY)_\a=\smatrix{0&0\\ Y_{\bar\a}&0}$.

We call a representation $X$  \emph{separated} if $(\Rad
X)_i=X_i$ for every sink $i$.

\begin{lem}
\pushQED{\qed}
Each representation $X$ admits a unique decomposition
\[X=X'\oplus\big(\bigoplus_{i\text{ a sink}}X(i)\big)\] such that $X'$ is
separated and each $X(i)$ is a direct sum of copies of $S(i)$.
\qedhere
\end{lem}

\begin{prop}\label{pr:radsq}
The functors $S$ and $T$ induce mutually inverse bijections between
the isomorphism classes of radical square zero representations of $Q$
and the isomorphism classes of separated representations of $Q^s$.
\end{prop}
\begin{proof}
It is easily checked that $TSX\cong X$ for each representation $X$ of
$Q$ satisfying $\Rad^2X=0$, and that $STY\cong Y$ for each separated
representation $Y$ of $Q^s$.
\end{proof}

\section{Representations of the Klein four group}

In this section we introduce representations of groups and sketch the
parallel with representations of quivers. In particular, we deduce the
classification of the representations for the Klein four group from
that for the Kronecker quiver.

\subsection{Representations of groups}

Let $G$ be a group and $k$ be a field. A \emph{$k$-linear
  representation} of $G$ is a pair consisting of a vector space $X$
over $k$ and a group homomorphism $G\to\Aut(X)$ into the group of
$k$-linear automorphisms of $X$, sending $g\in G$ to $X_g$.  One often
writes $gx=X_g(x)$ for $x\in X$.  A \emph{morphism} between two
representations $X,Y$ of $G$ is a $k$-linear map $\p\colon X\to Y$
such that $\p X_g=Y_g\p$ for all $g\in G$. We write $\Hom(X,Y)$ for
the set of all morphisms $X\to Y$.  The $k$-linear representations of
$G$ form a category which we denote by $\Rep(G,k)$.

Denote by $\Ga(G)$ the multiple loop quiver with $\Ga(G)_0=\{\ast\}$
and $\Ga(G)_1=G$.  Then the representations of $G$ can be identified
with representations $X$ of $\Ga(G)$ such
that \[X_1=\id_X\quad\text{and}\quad X_gX_h=X_{gh} \quad\text{for
  all}\quad g,h\in G.\] Thus all concepts developed for
representations of quivers carry over to representations of
groups.

\begin{exm} (1) The \emph{trival representation} $k$ consists of the
vector space $k$ with $gx=x$ for all $g\in G$ and $x\in k$.

(2) The \emph{regular representation} $k[G]$ consists of the vector
space $k[G]$ with basis $G$ and $g(\sum_h\a_hh)=\sum_h\a_hgh$ for all
$g\in G$ and $\sum_h\a_hh\in k[G]$
\end{exm}

\begin{lem}\label{le:Gproj}
Let $X$ be a representation of $G$. The map $\Hom(k[G],X)\to X$
sending $\p$ to $\p(1)$ is bijective.
\end{lem}
\begin{proof}
The inverse map sends $x\in X$ to $\p_x$ with $\p_x(\sum_g\a_g g)=\sum_g\a_g gx$
\end{proof}

The vector space duality $D=\Hom(-,k)$ induces a \emph{duality}
\[\Rep(G,k)\lto\Rep(G,k),\]
where $(DX)_g=D(X_{g^{-1}})$ for each $g\in G$; this yields an isomorphism
\[\Hom(X,DY)\cong\Hom(Y,DX)\] 
for each pair of representations $X,Y$.

\begin{lem}\label{le:Ginj}
Let $G$ be a finite group. Then the map $k[G]\to Dk[G]$ sending
$\sum_g\a_gg$ to $\sum_g\a_gg^*$ (with $g^*(h)=\d_{gh}$ for $g,h\in
G$) is an isomorphism.
\end{lem}
\begin{proof}
The map is bijective since $G$ is finite, and it is easily checked
that the map is a morphism, that is,  $\sum_g\a_g(hg)^*= h(\sum_g\a_gg^*)$ for all $h\in G$.
\end{proof}

\begin{rem}\label{re:Ginj}
It follows from Lemmas~\ref{le:Gproj} and \ref{le:Ginj} that the regular
representation $k[G]$ is a projective and injective object in the sense of
Remark~\ref{re:inj}.
\end{rem}

\subsection{Maschke's theorem}

The representation theory of a group depends heavily on the
characteristic of the field. A first instance of this phenomenon is
the following.

\begin{thm}[Maschke]
Let $X$ be a representation of a finite group $G$ such that the
characteristic of $k$ does not divide the order of $G$. Then every
subrepresentation $U\subseteq X$ admits a complement, that is, a
subrepresentation $V\subseteq X$ such that $X=U\oplus V$.  Therefore $X$
decomposes into a direct sum of simple representations.
\end{thm}
\begin{proof}
Let $r=\card G$. Choose a $k$-linear projection $\pi\colon X\to U$
onto $U$ and let $\pi'=r^{-1}\sum_{g\in G} X_g\pi X_{g^{-1}}$. Then
$\pi'$ is a morphism $X\to U$ such that $\pi'(x)=x$ for all
$x\in U$. Thus $X=U\oplus\Ker\pi'$.

The second assertion follows immediately by  induction on the
dimension of $X$.
\end{proof}

\subsection{Elementary abelian $p$-groups}
Let $p>0$ be prime and denote by $C_p$ the cyclic group of order
$p$.  For $r\ge 1$, we consider the \emph{elementary abelian group}
\[C_p^r=\langle g_1\ldots,g_r\rangle\]
and fix a field $k$ of characteristic $p$.  The representations of
$C_p^r$ can be identified with representations $X$ of the $r$-loop
quiver $\Ga$ (that is, $\Ga_0=\{\ast\}$ and
$\Ga_1=\{\g_1,\ldots,\g_r\}$) such that
\begin{equation}\label{eq:elabel}
X_{\g_i}X_{\g_j}=X_{\g_j}X_{\g_i}\quad\text{and}\quad X_{\g_i}^p=0
\quad\text{for all}\quad 1\le i,j\le r;
\end{equation}
the underlying vector space remains unchanged and
$X_{\g_i}=X_{g_i}-\id_X$ for each $i$.  This identification will be
used without further mentioning.

Recall that the Jacobson radical $\rad X$ of a representation $X$ is
the intersection of all its maximal subrepresentations.

\begin{lem}\label{le:rad_elabel}
Let $X$ be a representation of $C^r_p$. Then $\rad
X=\sum_{i=1}^r\Im X_{\g_i}$. Therefore $\rad^2 X=0$ if and only if
$X_{\g_i}X_{\g_j}=0$ for all $1\le i,j\le r$.  
\end{lem}
\begin{proof}
The relations in \eqref{eq:elabel} imply $\Rad^{pr}X=0$. Then it
follows from Lemma~\ref{le:jrad} that $\rad X=\Rad X=\sum_{i=1}^r\Im
X_{\g_i}$.
\end{proof}

We define a functor $T\colon\Rep(K_r,k)\to\Rep(C_p^r,k)\hookrightarrow\Rep(\Ga,k)$ from
representations of the $r$-Kronecker quiver 
\[\xymatrix{K_r&1\ar @/u2.75ex/@<3.0pt>[rr]^{\a_1}
\ar @/d2.75ex/@<-3.0pt>[rr]^{\vdots}_{\a_r}\ar[rr]^{\a_2}&& 2}\] to
representations of $C_p^r$. Given a Kronecker representation $X$, let
$TX=X_1\oplus X_2$ and $(TX)_{\g_i}=\smatrix{0&0\\ X_{\a_i}&0}$ for
$1\le i\le r$. For a morphism $\p\colon X\to Y$, let
$T\p=\smatrix{\p_1&0\\ 0&\p_2}$.

\begin{prop}\label{pr:elabel}
The functor $T$ induces a bijection between the isomorphism classes of
representations of the $r$-Kronecker quiver that have no direct
summand isomorphic to the simple representation $S(2)$, and the
isomorphism classes of representations $X$ of $C^r_p$ such that
$\rad^2 X=0$.
\end{prop}
\begin{proof}
The assertion is a special case of Proposition~\ref{pr:radsq} since
the $r$-Kronecker quiver is the separated quiver of $\Ga$. Note that a
Kronecker representation is separated if and only if there is no
direct summand isomorphic to $S(2)$. On the other hand, a
$\Ga$-representation $X$ satisfies $\Rad^2X=0$ if and only it is a
representation of $C^r_p$ satisfying $\rad^2 X=0$, by
Lemma~\ref{le:rad_elabel}.

Recall that the inverse of $T$ is constructed as follows. Let $X$ be a
representation of $C_p^r$ such that $X_{\g_i}X_{\g_j}=0$ for all $1\le
i,j\le r$. Then the $X_{\g_i}$ induce a family of linear maps
\[\xymatrix{X/\sum_i\Im X_{\g_i}\;\ar @/u0ex/@<8pt>^-{\bar X_{\g_1}}[rr]
\ar @/d0ex/@<-8pt>_-{\bar X_{\g_r}}[rr]^-{\vdots} &&\;\sum_i\Im
X_{\g_i}}.\] We denote this $r$-Kronecker quiver representation by
$SX$. It is easily checked that $TSX\cong X$, and that $STY\cong Y$
for each separated Kronecker representation $Y$.
\end{proof}

In view of Proposition~\ref{pr:wild_K3}, the category $\Rep(C_p^r,k)$
shows some wild behaviour for $r\ge 3$. More precisely, for any quiver
$Q$, there exists a functor \[F\colon\Rep(Q,k)\lto\Rep(C_p^r,k)\] such that
for every pair of representations $X,Y$ of $Q$,
\begin{enumerate}
\item $FX\cong FY$ implies $X\cong Y$, and
\item $F$ induces a monomorphism $\Hom(X,Y)\to\Hom(FX,FY)$.
\end{enumerate}
Note that one cannot expect a surjective map
$\Hom(X,Y)\to\Hom(FX,FY)$, since $\End(Z)\cong k$ implies $Z\cong k$ for each
representation $Z$ of $C_p^r$.

\subsection{The Klein four group}

Let $G=C_2\times C_2$ be the \emph{Klein four group}
and fix a field $k$ of characteristic $2$.

\begin{lem}\label{le:klein}
Let $X$ be a representation of $G$.  Then $X$ has no direct summand
isomorphic to $k[G]$ if and only if $X_{\g_1}X_{\g_2}=0$.
\end{lem}
\begin{proof}
Observe first that $1$, $\bar g_1=(g_1-1)$, $\bar g_2=(g_2-1)$, and
$\bar g_1\bar g_2=\bar g_2\bar g_1$ form a basis over $k$ of $k[G]$.  It follows that each
non-zero subrepresentation of $k[G]$ contains $\bar g_1\bar g_2$. 

Suppose there is an element $x\in X$ such that
$X_{\g_1}X_{\g_2}(x)\neq 0$. Using Lemma~\ref{le:Gproj}, there exists
a morphism $\p\colon k[G]\to X$ such that $\p(1)=x$. Then $\p(\bar g_1\bar g_2)=
X_{\g_1}X_{\g_2}(x)\neq 0$, and therefore $\Ker\p=0$. Thus $\p$ is a
split monomorphism, since $k[G]$ is an injective object; see
Remark~\ref{re:Ginj}. The other implication is clear.
\end{proof}

It follows from Proposition~\ref{pr:elabel} that the representations
of $G$ are closely related to representations of the Kronecker quiver
via the functor
$T\colon\Rep(K,k)\to\Rep(G,k)\hookrightarrow\Rep(\Ga,k)$, where
\[\xymatrix{K\qquad 1\ar@<3.0pt>^-{\a_1}[rr]\ar@<-3.0pt>_-{\a_2}[rr] &&2}
\qquad\text{and}\qquad
\xymatrix{\Ga\quad&\circ\ar@(ur,dr)^{\g_2}\ar@(ul,dl)_{\g_1}}.\] Given
Kronecker representations $X,Y$, let $TX=X_1\oplus X_2$ and
$(TX)_{\g_i}=\smatrix{0&0\\ X_{\a_i}&0}$ for $i=1,2$. For a morphism
$\p\colon X\to Y$, let $T\p=\smatrix{\p_1&0\\ 0&\p_2}$.

\begin{prop}\label{pr:klein}
The functor $T$ induces a bijection between the isomorphism classes of
representations of the Kronecker quiver that have no direct summand
isomorphic to the simple representation $S(2)$, and the isomorphism
classes of representations of the Klein four group that have no direct
summand isomorphic to the regular representation $k[G]$.
\end{prop}
\begin{proof}
The assertion is a special case of Proposition~\ref{pr:elabel}.  Note
that a representation $X$ of $G$ satisfies $\rad^2X=0$ if and only if
there is no direct summand isomorphic to $k[G]$, by
Lemmas~\ref{le:rad_elabel} and \ref{le:klein}.
\end{proof}

In view of Theorem~\ref{th:Kron}, we have the following
classification, assuming that the field $k$ is algebraically closed.

\begin{cor}
\pushQED{\qed} The representations $k[G]$, $TP_r$ ($r\ge 1$), $TI_r$
($r\ge 0$), and $TR_{p,\la}$ ($p\ge 1$ and $\la\in\bbP^1(k)$) form, up
to isomorphism, a complete list of finite dimensional indecomposable
representations of the Klein four group.\qedhere
\end{cor}

\section{Notes}

In \cite{G1} Gabriel proved that a connected quiver admits only
finitely many indecomposable representations if and only if it is of
Dynkin type. This result and the corresponding classification of the
indecomposable representations in terms of their dimension vectors is
known as Gabriel's theorem; see also \cite{G2}. The proof using
reflection functors is due to Bern{\v{s}}te\u{\i}n, Gel'fand, and
Ponomarev \cite{BGP}.

For the quivers of Euclidean type, the complete classification of
their indecomposable representations was established  by
Donovan and Freislich \cite{DF}, and independently by Nazarova
\cite{Na}; see also \cite{DR}. For the special case of the Kronecker
quiver, the classification is due to Kronecker \cite{Kr}.

The wild behaviour of representations of the $n$-subspace quiver for
$n\ge 5$ was first noticed by Brenner \cite{Br}.

The classification of the representations of the Klein four group was
obtained by Ba\v sev \cite{Ba}, and independently by Heller and Reiner
\cite{HR}. Gabriel's survey \cite{G2} discusses radical square zero
representations and further connections between representations of
finite groups and quivers.

Irreducible morphisms were introduced by Auslander and Reiten in
\cite{AR} and their structure for algebras of finite representation
type is discussed in \cite{G3}.  Infinite chains of morphisms were
used by Auslander when he characterised algebras of finite
representation type \cite{A}.

The exposition given here is based on various sources. Besides the
somewhat classical references given above, Crawley-Boevey's excellent
notes on representations of quivers \cite{CB} should be mentioned.

\begin{appendix}
\section{Exact sequences}

A sequence of morphisms between representations 
$$X_1\stackrel{\p_1}\lto
X_2\stackrel{\p_2}\lto\cdots\stackrel{\p_r}\lto X_{r+1}$$ is called
\emph{exact} if $\Im\p_i=\Ker\p_{i+1}$ for all $1\le i< r$.


\begin{Lem}\label{le:hom-seq}
\pushQED{\qed}
The following  are equivalent for a sequence $Y'\xto{\p}Y\xto{\psi}Y''$.
\begin{enumerate}
\item The sequence  $0\to Y'\xto{\p}Y\xto{\psi}Y''$ is exact.
\item The induced sequence  
$0\to \Hom(X,Y')\xto{(X,\p)}\Hom(X,Y)\xto{(X,\psi)}\Hom(X,Y'')$ is exact for each $X$.
\item The morphism $\p$ is a kernel of $\psi$.\qedhere
\end{enumerate}
\end{Lem}

\begin{Lem}
\pushQED{\qed}
The following  are equivalent for a sequence $Y'\xto{\p}Y\xto{\psi}Y''$.
\begin{enumerate}
\item The sequence  $Y'\xto{\p}Y\xto{\psi}Y''\to 0$ is exact.
\item The induced sequence  
$0\to \Hom(Y'',Z)\xto{(\psi,Z)}\Hom(Y,Z)\xto{(\p,Z)}\Hom(Y',Z)$ is exact for each $Z$.
\item The morphism $\psi$ is a cokernel of $\p$.\qedhere
\end{enumerate}
\end{Lem}

\begin{Lem}[Snake lemma]\label{le:snake}
Any commutative diagram 
\[\xymatrix{ 0\ar[r]&X'\ar[r]\ar[d]_{\p'}&X\ar[r]\ar[d]_{\p}&X''\ar[r]\ar[d]_{\p''}&0\\
 0\ar[r]&Y'\ar[r]&Y\ar[r]&Y''\ar[r]&0}\]
with exact rows induces an exact sequence
\[0\to\Ker\p'\to\Ker\p\to\Ker\p''\to\Coker\p'\to\Coker\p\to\Coker\p''\to 0.\]
\end{Lem}
\begin{proof}
The proof is given by Jill Clayburgh at the very beginning of the 1980
film \emph{It's My Turn} directed by Claudia Weill.
\end{proof}

An exact sequence 
\[0\lto X\stackrel{\p}\lto Y\stackrel{\psi}\lto Z\lto 0\]
is called \emph{split exact} provided that the following two
equivalent conditions are satisfied:
\begin{enumerate}
\item $\p$ is a split monomorphism, that is, $\id_X=\p'\p$ for some $\p'\colon Y\to X$;
\item $\psi$ is a split epimorphism, that is, $\id_Z=\psi\psi'$ for some $\psi'\colon Z\to Y$.
\end{enumerate}
We write  $\Ext(Z,X)=0$ if any exact sequence $0\to X\to Y\to Z\to 0$ splits.

\begin{Lem}\label{le:ext}
Let $0\to Y'\xto{}Y\xto{}Y''\to 0$ be an exact sequence and
$X$ a representation.
\begin{enumerate}
\item If $\Ext(Y',X)=0=\Ext(Y'',X)$, then $\Ext(Y,X)=0$.
\item If $\Ext(X,Y')=0=\Ext(X,Y'')$, then $\Ext(X,Y)=0$.
\end{enumerate}
\end{Lem}
\begin{proof}
We prove (1) and the argument for (2) is dual. Let $0\to X\xto{\p}
E\xto{\psi} Y\to 0$ be an exact sequence. Viewing the morphism $Y'\to Y$
as an inclusion gives rise to an exact sequence $0\to X\xto{}
E'\xto{} Y'\to 0$ by taking $E'=\psi^{-1}(Y')$. This sequence splits
and therefore $Y'\to Y$ factors through $\psi$. This yields the
following commutative diagram with exact rows
\[\xymatrix{&&Y'\ar[d]\ar@{=}[r]&Y'\ar[d]\\
0\ar[r]&X\ar[r]^\p\ar@{=}[d]&E\ar[r]^\psi\ar[d]&Y\ar[r]\ar[d]&0\\ 
0\ar[r]&X\ar[r]&E''\ar[r]&Y''\ar[r]&0}\]
by taking for $E''$ the cokernel of $Y'\to E$. The
lower sequence splits. Thus $X\to E''$ is a split monomorphism, and
therefore $\p$ is a split monomorphism. It follows that $\Ext(Y,X)=0$.
\end{proof}

\end{appendix}

\end{document}